\documentclass{amsart}

\usepackage{color}
\usepackage{url}

\newtheorem{theorem}{Theorem}[section]
\newtheorem{lemma}[theorem]{Lemma}

\theoremstyle{definition}

\theoremstyle{remark}
\newtheorem{remark}[theorem]{Remark}
\theoremstyle{corollary}
\newtheorem{corollary}[theorem]{Corollary}

\theoremstyle{conjecture}
\newtheorem{conjecture}[theorem]{Conjecture}
\numberwithin{equation}{section}

\def\k{\mathop{\mbox{\bf\large K}}}
\def\K{\mathop{\mbox{\bf\Large K}}}

\newcommand{\ba}{\boldsymbol{a}}

\begin{document}
\title{Continued fraction formulae involving ratios of three gamma functions}


\author{Xiaodong Cao}
\address{Department of Mathematics,
Beijing Institute of Petro-Chemical Technology,
Beijing, 102617, P.R. China}
\curraddr{}
\email{caoxiaodong@bipt.edu.cn}
\thanks{}

\author{Yoshio Tanigawa}
\address{Nishisato 2-13-1, Meitou, Nagoya 465-0084, Japan}
\curraddr{}
\email{tanigawa@math.nagoya-u.ac.jp,\ \ tanigawa${}_{-}$yoshio@yahoo.co.jp}
\thanks{}

\author{Wenguang Zhai}
\address{Department of Mathematics,
China University of Mining and Technology,
Beijing 100083, P.R. China}
\curraddr{}
\email{zhaiwg@hotmail.com}
\thanks{This work is supported by
the National Natural Science Foundation of China (Grant No.11971476)}

\subjclass[2010]{Primary 30B70; Secondary  33B15 41A25 68W25 11A55}

\date{}

\dedicatory{}
\keywords{Continued fraction; Gamma function; Bauer-Muir transformation; MC-algorithm; Modifying factors.}
\begin{abstract}
Via the MC-algorithm, in this paper we produce seven continued fraction formulae involving products and quotients of three gamma functions with three parameters, and another is an extension of Entry 34 in Chapter 12 of Ramanujan's second notebook. Five of them will be proved rigorously by the Bauer-Muir transformation. A crucial ingredient in the proofs of our five theorems is to employ the Bauer-Muir transformation twice with different nonlinear modifying factors.
\end{abstract}

\maketitle
\section{Introduction and results}
Continued fractions belong to classical research areas of mathematics. An intensive study of continued fractions started since John Wallis and Lord Brouncker. Many mathematicians such as Cauchy, Euler, Gauss, Hermite, Jacobi, Perron, Stieltjes, etc. investigated continued fractions~\cite{Brezinski}. Ramanujan~(see~\cite{ABBW,ABJL,Berndt,BLW,Ram}) produced quite a number of continued fraction expansions for products and quotients of two, four, or eight gamma functions. These formulae have all proved to be connected to hypergeometric functions, for example, see \cite[p.~284]{LW1}, \cite{Zhang}. In Ramanujan's notebooks, however, three of principle formulae involving gamma functions are Entries 34, 39, and 40 in Chapter 12 of Ramanujan's second notebook, for instance, see Berndt et al.~\cite{BLW}, Masson~\cite{Masson,Masson1991} and Watson~\cite{Watson}. Indeed, Ramanujan's contributions to the continued fraction expansions of analytic functions are one of his most spectacular achievements. Unfortunately, Ramanujan left us no clues as to how he discovered these elegant continued fraction formulae~(see \cite{Berndt, BLW}). 
For more details of the history of Ramanujan's continued fraction formulae for products and quotients of gamma functions, interested reader is referred to the Ph.D. thesis of Reuter~\cite{Reuter} and references quoted therein. It is worthy to note that the recent pioneering work, dubbed ``The Ramanujan Machine"  by Raayoni et al.~\cite{Raayoni,Raayoni-1}, produces many interesting continued fraction formulae, some of which are difficult to prove true. Dougherty-Bliss and Zeilberger~\cite{DZ} develop their beautiful approach, they use \emph{symbolic} rather than \emph{numeric} computation. However, their approaches are useful only for at most one parameter. The main purpose of this paper is to establish five continued fraction expansions for products and quotients of three gamma functions with three parameters. Our results will provide some new structures or models for discovering continued fraction formulae in the future.

To describe our results more clearly, let us prepare some notation. We use $\mathbb{N}_0$ to denote the set of non-negative integers. The set of integers is denoted by $\mathbb{Z}$. The notation $P_l(x)$ or $Q_l(x)$ denote a polynomial of degree $l$ in $x$. We shall employ $[x^j]P(x)$ to denote the coefficient of $x^j$ in polynomial $P(x)$. Assume that $\{a_n\}_{n=1}^\infty$ and $\{b_n\}_{n=0}^\infty$ are two sequences of complex numbers. The generalized continued fraction
\begin{align}
\tau=b_0+\frac{a_1}{b_1+\frac{a_2}{b_2+\ddots}}=b_0+
\begin{array}{ccccc}
a_1 && a_2 &       \\
\cline{1-1}\cline{3-3}\cline{5-5}
 b_1 &\!+\!& b_2 &\!+\cdots
\end{array}
=b_0+\K_{n=1}^{\infty}
\Bigl(\frac{a_n}{b_n}\Bigr)
\end{align}
is defined as the limit of the $n$th approximant
\begin{equation}
\frac{A_n}{B_n}=b_0+\K_{k=1}^{n}\Bigl(\frac{a_k}{b_k}\Bigr)
\end{equation}
as $n$ tends to infinity. We shall employ the space-saving notation
\begin{equation}
\frac{a_0}{b_0+}\K_{m=1}^{\infty}\Bigl(\frac{a_m}{b_m}\Bigr):=\frac{a_0}{b_0+\K_{m=1}^{\infty}\Bigl(\frac{a_m}{b_m}\Bigr)}.
\end{equation}
We refer to the books by Cuyt et al.~\cite{CPV}, Jones and Thron~\cite{JT}, Khovanskii~\cite{Kh}, Lorentzen and Waadeland~\cite{LW,LW1}, Perron~\cite{Perron}, Wall~\cite{Wall}, for the present state of the theory of continued fractions. 
The main results of this paper are stated as follows.

\begin{theorem}
Assume that $x$ is complex with $\Re(x)>0$, or assume that either one of $\alpha, \beta$ is an odd integer or $\alpha+\beta$ is an even integer. Then for $j=1,2,3$ 
\begin{align}
&\frac{ \Gamma\left(\frac 12 ( x+\alpha)\right) \Gamma\left(\frac 12(x-\alpha-\beta+1)\right)\Gamma
  \left(\frac 12 (x +\beta)\right)}{
\Gamma\left(\frac 12 (x-\alpha+2)\right) \Gamma\left(\frac 12( x+\alpha+\beta+1)\right)\Gamma
 \left( \frac 12 (x-\beta+2)\right)}\label{Theorem 1-First}\\
=&\frac{2}{\lambda_0(j;x)/2+}\K_{m=1}^{\infty}\Bigl(\frac{\kappa_m(j)}{\lambda_m(j;x)}\Bigr)
=\frac{2}{(x^2+q_0)/2+}\K_{m=1}^{\infty}\Bigl(\frac{p_m}{x^2+q_m}\Bigr),\label{Theorem 1-Second}
\end{align}
where  for $m\in \mathbb{N}_0$
\begin{equation}\begin{cases}
\kappa_{2m}(1)=\frac{(2m-2+\alpha+\beta)(2m-1-\alpha)(2m-1-\beta)}{4(2m-1)},\\
\kappa_{2m-1}(1)=\frac{(2m-\alpha-\beta)(2m-1+\alpha)(2m-1+\beta)}{4(2m-1)},
\end{cases}
\begin{cases}
\lambda_{2m}(1;x)=x^2-1,\\
\lambda_{2m-1}(1;x)=1,
\end{cases}
\end{equation}
\begin{align}
&\begin{cases}
\kappa_{2m}(2)=\frac{(2m-2+\alpha+\beta)(2m-1+\alpha)(2m-1+\beta)}{4(2m-1)},\\
\kappa_{2m-1}(2)=\frac{(2m-\alpha-\beta)(2m-1-\alpha)(2m-1-\beta)}{4(2m-1)},
\end{cases}\\
&\begin{cases}
\lambda_{2m}(2;x)=x^2-(1-\alpha-\beta)^2,\\
\lambda_{2m-1}(2;x)=1;
\end{cases}
\end{align}
\begin{align}
\begin{cases}
\kappa_{2m}(3)=\frac{(2m-\alpha-\beta)(2m-1+\alpha)(2m-1-\beta)}{4(2m-1)},\\
\kappa_{2m-1}(3)=\frac{(2m-2+\alpha+\beta)(2m-1-\alpha)(2m-1+\beta)}{4(2m-1)},
\end{cases}
\begin{cases}
\lambda_{2m}(3;x)=x^2-\beta^2,\\
\lambda_{2m-1}(3;x)=1;
\end{cases}
\end{align}
and
\begin{align}\begin{cases}
p_m=-\frac{(2m-\alpha-\beta)(2m-2+\alpha+\beta)\left((2m-1)^2-\alpha^2\right)\left((2m-1)^2-\beta^2\right)}{16(2m-1)^2},\\
q_m=\frac 12\left(4m^2+\alpha-\alpha^2+\beta-\beta^2-\alpha\beta\right)+\frac{\alpha\beta(-1+\alpha+\beta)}{2(2m-1)(2m+1)}.
\end{cases}\end{align}
\end{theorem}

\begin{theorem}
Suppose that $x$ is complex with $\Re(x)>0$, or suppose that either one of $\alpha, \beta$ is an odd integer or $\alpha+\beta$ is an even integer. We then have
\begin{align}
&\frac{ \Gamma\left(\frac 12 ( x - \alpha+1)\right) \Gamma\left(\frac 12(x+\alpha+\beta)\right)\Gamma
  \left(\frac 12 (x -\beta+1)\right)}{
\Gamma\left(\frac 12 ( x +\alpha+1)\right) \Gamma\left(\frac 12( x-\alpha-\beta+2)\right)\Gamma
 \left( \frac 12 (x +\beta+1)\right)}\label{Theorem 2}\\
=&\frac{1}{x/2+}\K_{m=1}^{\infty}\Bigl(\frac{\kappa_m}{x}\Bigr),\nonumber
\end{align}
where for $m\in \mathbb{N}$
\begin{equation}\begin{cases}
\kappa_{2m}=\frac{(2m-\alpha-\beta)(2m-1+\alpha)(2m-1+\beta)}{4(2m-1)},\\
\kappa_{2m-1}=\frac{(2m-2+\alpha+\beta)(2m-1-\alpha)(2m-1-\beta)}{4(2m-1)}.
\end{cases}
\end{equation}
\end{theorem}

\begin{theorem}
Assume that either $x$ is complex with $\Re(x)>\Re\left(1-\alpha-\beta\right)$, or one of $\alpha, \beta, \alpha+\beta$ is an even integer. Then
\begin{align}
&\frac{ \Gamma\left(\frac 12 ( x - 2\alpha-\beta+2)\right) \Gamma\left(\frac 12(x+2)\right)\Gamma
  \left(\frac 12 (x -\alpha-2\beta+2)\right)}{
\Gamma\left(\frac 12 ( x-\alpha+2)\right) \Gamma\left(\frac 12(x-2(\alpha+\beta)+4)\right)\Gamma
 \left( \frac 12 (x-\beta+2)\right)}\label{Theorem 3}\\
=&\frac{1}{(x+\lambda_0)/2+}\K_{m=1}^{\infty}\Bigl(\frac{\kappa_m}{x+\lambda_m}\Bigr),\nonumber
\end{align}
where for $m\in \mathbb{N}_0$
\begin{equation}\begin{cases}
\kappa_{2m}=\frac{(2m-\alpha-\beta)(2m-\alpha)(2m-\beta)}{4(2m-1)},\\
\kappa_{2m-1}=\frac{(2m-2+\alpha+\beta)(2m-2+\alpha)(2m-2+\beta)}{4(2m-1)},
\end{cases}\quad
\begin{cases}
\lambda_{2m}=2(1-\alpha-\beta),\\
\lambda_{2m-1}=0.
\end{cases}
\end{equation}
\end{theorem}

\begin{theorem}
Assume that either $x$ is complex with $\Re(x)>\Re\left(1-\alpha-\beta\right)$, or one of $\alpha, \beta, \alpha+\beta$ is an even integer. Then
\begin{align}
&\frac{ \Gamma\left(\frac 12 ( x - 2\alpha-\beta+2)\right) \Gamma\left(\frac 12x\right)\Gamma
  \left(\frac 12 (x -\alpha-2\beta+2)\right)}{
\Gamma\left(\frac 12 ( x - \alpha+2)\right) \Gamma\left(\frac 12( x-2(\alpha+\beta)+2)\right)\Gamma
 \left( \frac 12 (x -\beta+2)\right)}\label{Theorem 4}\\
=&\frac{1}{(x+\lambda_0)/2+}\K_{m=1}^{\infty}\!\Bigl(\frac{\kappa_m}{x+\lambda_m}\Bigr),\nonumber
\end{align}
where for $m\in \mathbb{N}_0$
\begin{equation}\begin{cases}
\kappa_{2m}=\frac{(2m-\alpha-\beta)(2m-\alpha)(2m-\beta)}{4(2m-1)},\\
\kappa_{2m-1}=\frac{(2m-2+\alpha+\beta)(2m-2+\alpha)(2m-2+\beta)}{4(2m-1)},
\end{cases}
\begin{cases}
\lambda_{2m}=0,\\
\lambda_{2m-1}=2(1-\alpha-\beta).
\end{cases}
\end{equation}
\end{theorem}


\begin{theorem}
Let $x, \alpha, \beta\in \mathbb{C}$, $\ba=(\alpha,\beta)$, and $P(\ba;x)$ be defined by
\begin{equation}
P(\ba;x):=\frac{ \Gamma\left(\frac 12 ( x - \alpha+1)\right) \Gamma\left(\frac 12(x+\alpha+\beta+1)\right)\Gamma
  \left(\frac 12 (x -\beta+1)\right)}{
\Gamma\left(\frac 12 ( x +\alpha+1)\right) \Gamma\left(\frac 12( x-\alpha-\beta+1)\right)\Gamma
 \left( \frac 12 (x +\beta+1)\right)}.\label{definition of P}
\end{equation}
Assume that $x$ is complex such that $\Re(x)>0$, or assume that either $\alpha$, $\beta$, or $\alpha+\beta$ is an even integer. Then
\begin{align}
\frac{1-P(\ba;x)}{1+P(\ba;x)}=&\frac{\alpha\beta(\alpha+\beta)/4}{b_0(x)+}\K_{m=1}^{\infty}\Bigl(\frac{a_m}{b_m(x)}\Bigr)\label{first assertion of theorem 5}\\
=&\frac{\alpha\beta(\alpha+\beta)/4}{x^2+\lambda_0+}\K_{m=1}^{\infty}\Bigl(\frac{\kappa_m}{x^2+\lambda_m}\Bigr),\label{second assertion of theorem 5}
\end{align}
where $b_0(x)=x^2-1-\alpha\beta(\alpha+\beta)/4$, and for $k\in\mathbb{N}$
\begin{align}
&\begin{cases}
a_{2m}=\frac{1}{4}(2m+\alpha)(2m-\alpha-\beta)(2m+\beta),\\
a_{2m-1}=
\frac 14{(2m-\alpha)(2m+\alpha+\beta)(2m-\beta)},
\end{cases}\\
&\begin{cases}
b_{2m}(x)=(x^2-1)(2m+1),\\
b_{2m-1}(x)=1;
\end{cases}
\end{align}
\begin{equation}\begin{cases}
\kappa_m=-\frac{\left((2m)^2-\alpha^2\right)\left((2m)^2-(\alpha+\beta)^2\right)\left((2m)^2-\beta^2\right)}{16 (2 m - 1) (2 m +
     1)},\quad & m\ge 1,\\
\lambda_m=\frac 12 \left(4 m^2+4m+2 - \alpha^2-\alpha\beta-\beta^2\right), \quad & m\ge 0.
\end{cases}
\end{equation}
\end{theorem}

On taking $\alpha=\beta=1/3$ in Theorem 1, we deduce the following continued fraction formula.

\begin{corollary} For $m\in\mathbb{N}$, let $\kappa_{2m}=\frac{2(3m-2)(3m-1)^2}{27(2m-1)}$, $\kappa_{2m-1}=\frac{2(3m-1)(3m-2)^2}{27(2m-1)}$,
$\nu_{2m}=\frac{2(3m-2)^3}{27(2m-1)}$ and $\nu_{2m-1}=\frac{2(3m-1)^3}{27(2m-1)}$. If $\Re(x)>0$, then
\begin{align}
&\frac{\Gamma^3\big(\frac 12 (x+\frac 13)\big)}{\Gamma^3\big(\frac 12 (x+\frac 53)\big)}\label{Coro-1-1}\\
=&\begin{array}{ccccccccccc}
2& & \kappa_1& & \kappa_2&&\kappa_{2m-1}&&\kappa_{2m}&\\
\cline{1-1}\cline{3-3}\cline{5-5}\cline{7-7}\cline{9-9}\cline{11-11}
\big(x^2-1/9\big)/2&\!+\!&1&\!+\!&x^2-1/9&\!+\!\cdots\!+\!&1&\!+\!&x^2-1/9&\!+\!\cdots\!
\end{array}\nonumber\\
=&\begin{array}{ccccccccccc}
2& & \nu_1& & \nu_2&&\nu_{2m-1}&&\nu_{2m}&\\
\cline{1-1}\cline{3-3}\cline{5-5}\cline{7-7}\cline{9-9}\cline{11-11}
\big(x^2-1\big)/2&\!+\!&1&\!+\!&x^2-1&\!+\!\cdots\!+\!&1&\!+\!&x^2-1&\!+\!\cdots\!
\end{array}\nonumber
\end{align}
\end{corollary}
Similarly, putting $\alpha=\beta=1/3$ in Theorem 2, we obtain
\begin{corollary} For $m\in\mathbb{N}$, we let $\kappa_{2m}=\frac{2(3m-1)^3}{27(2m-1)}$ and $\kappa_{2m-1}=\frac{2(3m-2)^3}{27(2m-1)}$. If $\Re(x)>0$, then
\begin{equation}
\frac{\Gamma^3\big(\frac 12 (x+\frac 23)\big)}{\Gamma^3\big(\frac 12 (x+\frac 43)\big)}=\frac{1}{x/2+}\K_{m=1}^{\infty}\Bigl(\frac{\kappa_m}{x}\Bigr).\label{Coro-2}
\end{equation}
\end{corollary}
\begin{remark}
For $n\in\mathbb{N}$, let
\begin{align}\begin{cases}
T_1(n):=\frac{1\cdot 4\cdots (3n-2)}{3\cdot6\cdots (3n)}
=\frac{\Gamma(n+\frac 13)}{\Gamma(n+1)\Gamma(\frac 13)},\\
T_2(n):=\frac{2\cdot 5\cdots (3n-1)}{3\cdot6\cdots (3n)}=\frac{\Gamma(n+\frac 23)}{\Gamma(n+1)\Gamma(\frac 23)}.
\end{cases}
\end{align}
See \cite[p.~255, (6.1.11) and (6.1.13)]{AS}. Mortici et al.~\cite{MCL} first investigated the continued fraction approximations for $T_1$ and $T_2$. Bai et al.~\cite{BLW} showed some inequalities for $T_1$ and $T_2$. Corollary 1 and 2 confirm two questions proposed in \cite{CTZ,CW}.
We may compare the above two corollaries with the following Bauer's continued fraction in 1872, which was communicated by Ramanujan~\cite[p.~xxvii]{Ram62} in his first letter to Hardy:
\begin{align}
\frac{\Gamma^2\big(\frac 14 (x+1)\big)}{\Gamma^2\big(\frac 14(x+3)\big)}
=\begin{array}{ccccccccccc}
4& & 1^2& & 3^2&&5^2&&7^2&\\
\cline{1-1}\cline{3-3}\cline{5-5}\cline{7-7}\cline{9-9}\cline{11-11}
x&\!+\!&2x&\!+\!&2x&+&2x&\!+\!&2x&\!+\!\cdots\!,\quad \Re(x)>0.
\end{array}
\end{align}
\end{remark}

Replacing $\alpha$ and $\beta$ all by $\alpha+2/3$ in Theorem 3 and 4, respectively, we establish the following two formulae.
\begin{corollary}  Suppose that either $\Re(x)>-\Re(2\alpha+1/3)$ or one of $\alpha+2/3$ and $-\alpha+1/3$ is a positive integer. Then
\begin{align}
&\frac{\Gamma\big(\frac 12 (x-\alpha)\big)\Gamma\big(\frac 12 (x+2)\big)\Gamma\big(\frac 12 (x+\alpha)\big)}{\Gamma\big(\frac 12 (x-\alpha+\frac 43)\big)\Gamma\big(\frac 12 (x+\frac 43)\big)\Gamma\big(\frac 12 (x+\alpha+\frac 43)\big)}\label{Coro 3-First}\\
=&\frac{1}{(x+\lambda_0)/2+}\K_{m=1}^{\infty}\Bigl(\frac{\kappa_m}{x+\lambda_m}\Bigr),\nonumber\\%
&\frac{\Gamma\big(\frac 12 (x-\alpha)\big)\Gamma\big(\frac 12 x\big)\Gamma\big(\frac 12 (x+\alpha)\big)}{\Gamma\big(\frac 12 (x-\alpha+\frac 43)\big)\Gamma\big(\frac 12 (x-\frac 23)\big)\Gamma\big(\frac 12 (x+\alpha+\frac 43)\big)}\label{Coro 3-Second}\\
=&\frac{1}{(x+\mu_0)/2+}\K_{m=1}^{\infty}\Bigl(\frac{\kappa_m}{x+\mu_m}\Bigr),\nonumber
\end{align}
where
\begin{align}
&\begin{cases}
\kappa_{2m}=\frac{\left(m-\alpha-2/3\right)\left(2m-\alpha-2/3\right)^2}{2(2m-1)},\\
\kappa_{2m-1}=\frac{\left(m+\alpha-1/3\right)\left(2m+\alpha-4/3\right)^2}{2(2m-1)},
\end{cases}\\
&\begin{cases}
\lambda_{2m}=-2\left(2\alpha+\frac 13\right),\\
\lambda_{2m-1}=0,
\end{cases}
\begin{cases}
\mu_{2m}=0,\\
\mu_{2m-1}=-2\left(2\alpha+\frac 13\right).
\end{cases}
\end{align}
\end{corollary}
The paper is organized as follows. In section 2, based on the work on~\cite{CTZ}, we shall further develop the multiple-correction method for function approximation, or the MC-algorithm for short, which may be used to guess continued fraction formulae for a class of functions with multiple parameters. In section 3, we shall prepare the Bauer-Muir transformation in the theory of continued fractions, which plays an important role in the proofs of Theorem 1 to 5. From Sec. 4 to 8, we shall sequentially give the proofs of five theorems. In the last section, we shall propose three conjectures for further research.
\section{The MC-algorithm for function approximation}
A large part of this section is taken from~\cite{CTZ} which, in turn, follows from earlier works of Mortici~\cite{Mor} and Cao~\cite{Cao1}. Let $f(x)$ be a real function defined on $(x_0,+\infty)$ to be approximated with $\lim_{x\rightarrow+\infty}f(x)=0$.  Our main interest is to find a better rational approximation of $f(x)$ for those large $x$. Throughout this section, we assume that there exists a fixed positive integer $\nu$ and a constant $c\neq 0$ such that $\lim_{x\rightarrow +\infty}x^{\nu}f(x)=c$. In this case, we say that the function $f(x)$ is of order $x^{-\nu}$ when $x$ tends to infinity, and we employ the notation $\mathrm{R}(f(x))$ to denote the exponent $\nu$ of $x^{\nu}$. For convenience, $\mathrm{R}(0)$ is stipulated to be infinity. In a certain sense, $\mathrm{R}(f(x))$ characterizes the rate of convergence for $f(x)$ when $x$ tends to infinity. Inspired by Pad\'e approximation's idea~(~for example, see \cite[p.~35, \S 1.4.4]{LW1}~), the question mentioned above can be described as: given a fixed positive $k$, how to find two polynomials $P_l(x)$ and $Q_m(x)$ such that $\mathrm{R}\Bigl(f(x)-\frac{P_l(x)}{Q_m(x)}\Bigr), (l<m\le k)$  attains maximum, i.e.
\begin{equation}
\max\limits_{P_l(x),Q_m(x)}\left\{
\mathrm{R}\Bigl(f(x)-\frac{P_l(x)}{Q_m(x)}\Bigr): l<m\le k\right\}.
\end{equation}
Let $\Phi_j(\kappa_j;x)$ or $\Psi_j(\kappa_j;x)$ denote a monic polynomial of degree $\kappa_j$ in $x$. By means of Euclidean algorithm, however, it is not difficult to show that any rational function $P_l(x)/Q_m(x)$ with $l<m$ can be uniquely written in the form
\begin{equation}
\frac{P_l(x)}{Q_m(x)}=\K_{j=0}^J\Bigl(\frac{\lambda_j}{\Phi_j(\kappa_j;x)}\Bigr),
\end{equation}
where $\lambda_j\in \mathbb{R}\setminus\{0\}$. This is the main reason why we always consider the continued fraction approximation of $f(x)$. More precisely, for the MC-algorithm, instead of finding $P_l(x)/Q_m(x)$, we try to solve sequentially $\left(\lambda_0,\Phi_0(\kappa_0;x)\right)$,
$\left(\lambda_1,\Phi_1(\kappa_1;x)\right)$,\ldots, $\left(\lambda_{k^*},\Phi_{k^*}(\kappa_{k^*};x)\right)$ until some suitable $k^*$ you want.

We shall now describe how to solve the $k$th correction $\mathrm{MC}_k(x):=\k_{j=0}^k\Bigl(\frac{\lambda_j}{\Phi_j(\kappa_j;x)}\Bigr)$. Firstly, the initial-correction $\mathrm{MC}_0(x)=\frac{\lambda_0}{\Phi(\kappa_0;x)}$ is essential. If $f(x)\sim \frac{\lambda_0}{x^{\kappa_0}}$ as $x$ tends to infinity, we choose
$\Phi(\kappa_0;x)$ such that
\begin{equation}
\mathrm{R}\Bigl(f(x)-\frac{\lambda_0}{\Phi(\kappa_0;x)}\Bigr)=
\max\limits_{\Psi(\kappa_0; x)}\left\{
\mathrm{R}\Bigl(f(x)-\frac{\lambda_0}{\Psi(\kappa_0;x)}\Bigr)\right\}.
\end{equation}
Secondly, in order to solve the first-correction $\mathrm{MC}_1(x)=\frac{\lambda_0}{\Phi(\kappa_0;x)+}\frac{\lambda_1}{\Phi_1(\kappa_1;x)}$, first we find the smallest positive integer $\kappa_1$ and $\lambda_1\neq 0$ such that
\begin{equation}
\mathrm{R}\Bigl(f(x)-\frac{\lambda_0}{\Phi(\kappa_0;x)+\frac{\lambda_1}{x^{\kappa_1}}}\Bigr)<
\mathrm{R}\Bigl(f(x)-\frac{\lambda_0}{\Phi(\kappa_0;x)}\Bigr).
\end{equation}
Then, $\Phi_1(\kappa_1;x)$ can be established by the following relation
\begin{align}
&\mathrm{R}\Bigl(f(x)-\frac{\lambda_0}{\Phi(\kappa_0;x)+}\frac{\lambda_1}{\Phi_1(\kappa_1;x)}\Bigr)\\
=&
\max\limits_{\Psi(\kappa_1;x)}\left\{
\mathrm{R}\Bigl(f(x)-\frac{\lambda_0}{\Phi(\kappa_0;x)+}\frac{\lambda_1}{\Psi(\kappa_1;x)}\Bigr)\right\}.\nonumber
\end{align}
Finally, we may repeat the above process to sequentially find the $k$th correction $\mathrm{MC}_k(x)$. Roughly speaking, a key step of the MC-algorithm is to use a polynomial $\Psi(l;x)$ of degree $l$ instead of $x^{l}$ for acceleration of convergence.
\subsection{The Mortici lemma}
Mortici~\cite{Mor} established a very useful tool for measuring the rate of convergence, which claims that a sequence $\{x_n\}_{n=1}^{\infty}$ converging to zero is the fastest possible when the difference $\{x_n-x_{n+1}\}_{n=1}^{\infty}$ is the fastest possible. Since then, the Mortici lemma has been effectively applied in many papers such as~\cite{BDL,Cao1,Cao2,CTZ,CW,MCL}. In fact, the Mortici lemma may be viewed as a special form of Stolz formula. The following lemma is a simple generalization of the Mortici lemma. For the details of the proof , the reader may refer to~\cite{Cao2}.
\begin{lemma}(The Mortici lemma) If $\lim_{x\rightarrow+\infty}f(x)=0$, and there exists the limit
\begin{align}
 \lim_{x\rightarrow+\infty}x^\lambda\left(f(x)-f(x+1)\right)=l\in
 \mathbb{R},
\end{align}
with $\lambda>1$, then
\begin{align}
 \lim_{x\rightarrow+\infty}x^{\lambda-1}f(x)=\frac{l}{\lambda-1}.
\end{align}
\end{lemma}

\subsection{The Mortici transformation}
In this subsection, we shall explain how to find all the related coefficients in $\mathrm{MC}_k(x)$. If we can expand $f(x)$ into a power series in $1/x$ easily, it is not difficult to determine $\mathrm{MC}_k(x)$. Similarly, if we may expand the difference $f(x)-f(x+1)$ into a power series in $1/x$, by the Mortici lemma we can also find $\mathrm{MC}_k(x)$, e.g. the Euler-Mascheroni
constant, the constants of Landau, the constants of Lebesgue, etc.~(see~\cite{Cao1}). However, in some cases the previous two approaches are not very efficient, for instance, the ratios of the gamma functions. We shall employ the following approach developed by Mortici~\cite{Mor} to treat such a situation.

First we introduce the $k$th-correction relative error sequence $\{E_k(x)\}_{k=0}^{\infty}$ to be defined by the following recurrence relations
\begin{equation}
f(x)=\mathrm{MC}_0(x)
\exp\left(E_0(x)\right),\quad
f(x)=\mathrm{MC}_k(x)\exp\left(E_k(x)\right),~~(k\ge 1).\label{Ek-def}
\end{equation}
We see easily that for $k\in \mathbb{N}_0$,
\begin{equation}f(x)-\mathrm{MC}_k(x)=\mathrm{MC}_k(x)\left(\exp\left(E_k(x)\right)-1\right).
\end{equation}
It follows from $\lim_{x\rightarrow\infty}E_k(x)=0$ and $\lim_{t\rightarrow 0}\frac{\exp(t)-1}{t}=1$ that for $k\in \mathbb{N}_0$,
\begin{align}
\mathrm{R}\left(f(x)-\mathrm{MC}_k(x)\right)
=\kappa_0+\mathrm{R}\left(E_k(x)\right).\label{Mortici-relation-1}
\end{align}
In this way, we reduce the problem to solving $\mathrm{R}\left(E_k(x)\right)$.

Take the logarithm of two identities in~\eqref{Ek-def}, respectively, we deduce that  for $k\in \mathbb{N}_0$
\begin{equation}
\ln f(x) =\ln\mathrm{MC}_k(x)+E_k(x).
\end{equation}

Next, let us consider the difference
\begin{equation}
E_k(x)-E_k(x+1)=\ln\frac{f(x)}{f(x+1)}
+\ln\frac{\mathrm{MC}_k(x+1)}{\mathrm{MC}_k(x)},\quad k\in \mathbb{N}_0.\label{Mortici-relation-2}
\end{equation}
By means of the Mortici lemma, one gets
\begin{equation}
\mathrm{R}\left(E_k(x)\right)=\mathrm{R}\left(E_k(x)-E_k(x+1)\right)-1.\label{Mortici-relation-3}
\end{equation}

Combining \eqref{Mortici-relation-1}, \eqref{Mortici-relation-2},and \eqref{Mortici-relation-3}, we attain the following useful tool.
\begin{lemma}(The Mortici transformation)  With the above notation, for $k\in \mathbb{N}_0$
\begin{equation}
\mathrm{R}\left(f(x)-\mathrm{MC}_k(x)\right)
=-1+\kappa_0+\mathrm{R}\Bigl(\ln\frac{f(x)}{f(x+1)}
+\ln\frac{\mathrm{MC}_k(x+1)}{\mathrm{MC}_k(x)}\Bigr).\label{Mortici-transformation}
\end{equation}
\end{lemma}
\begin{remark} The idea of Lemma 2 is originated from Mortici~\cite{Mor}, which will be called a Mortici transformation. In its applications, it should be stressed that we always choose a suitable $f(x)$ such that $f(x)/f(x+1)$ is a rational function in $x$ by using the recurrence relation $\Gamma(x+1)=x\Gamma(x)$ and variable transformation.
\end{remark}
\begin{remark}
To guess continued fraction formula of the form $\frac{1-P}{1+P}$, we often use the following simple transformation~(see Lorentzen and Waadeland~\cite[p.~560, (1.1.2)]{LW})
\begin{equation}\label{Usual simple transformation}
\frac{1-P}{1+P}=1-\frac{2}{1+b_0+}\K_{m=1}^{\infty}\Bigl(\frac{a_m}{b_m}\Bigr)\Leftrightarrow P=b_0+\K_{m=1}^{\infty}\Bigl(\frac{a_m}{b_m}\Bigr).
\end{equation}
For instance, Entries 34 and 35 in Chapter 12 of Ramanujan's second notebook~\cite{Berndt}, Theorem 5 and three conjectures in Section 9 below.
With the aid of the MC-algorithm, we can guess almost all continued fraction formulae in Chapter 12 of Ramanujan's notebooks~(see~\cite{Berndt, BLW}) involving products and quotients of gamma functions apart from his Entry 40.
\end{remark}


On the one hand, in order to determine all related coefficients, we often use an appropriate symbolic computation software, which requires a huge amount of computations. On the other hand, the exact expression at each occurrence also takes a lot of space. All theorems in this  paper are built on experimental results in this section. Hence, we shall focus on the rigorous proofs of main results, and omit some of the related details for guessing these formulae. Interested readers may refer to Sec. 6 and 8 in preprint~\cite{CTZ}.
\section{The Bauer-Muir transformation}
One of the outstanding theorems in the theory of continued fractions is the result described by O. Perron~\cite{Perron} as the Bauer-Muir transformation. Although the Bauer-Muir transformation is too complicated and amazing, this transformation and limiting cases of it give rise to numerous and extremely interesting
consequences.  For instance, see Andrews et al.~\cite{ABJL}, Berndt et al.~\cite{BLW}, Berndt~\cite[Chapter 32]{Berndt98}, Cao and Chen~\cite{CC}, Jacobsen~\cite{Jac-1}, Lorentzen and Waadeland~\cite[p.~77--82, Example 11 to 13]{LW}, and Perron~\cite[p.~25-37]{Perron}, etc. We shall use the Jacobsen's formulation in ~\cite{Jac-1}~(also see~\cite[p.~76, Theorem 11]{LW} or \cite[p.~83, Theorem 2.18]{LW1}). Thanks to Lorentzen for a series of important work~\cite{Jac,Jac-1,Loren94c}, which plays a key role in proofs of our five theorems with non-positive continued fractions or non-positive modifying factors. 
\begin{lemma} (The Bauer-Muir transformation) The Bauer-Muir transformation of  $b_0+\k\big(a_m/b_m\big)$ with respect to
$\left\{r_m\right\}_{m=0}^{\infty}$ from $\mathbb{C}$ exists if and only if
\begin{align}
\varphi_m=a_m-r_{m-1}\left(b_m+r_m\right)\neq 0,\quad m\ge 1.
\end{align}
If it exists, then it is given by
\begin{align*}
&b_0+\begin{array}{ccccccc}
a_1& & a_2& & a_3 &\\
\cline{1-1}\cline{3-3}\cline{5-5}\cline{7-7}
b_1&+&b2&+&b_3&
+\cdots
\end{array}\\
=&b_0+r_0+
\begin{array}{ccccccc}
\varphi_1& & a_1 \varphi_2/\varphi_1& & a_2 \varphi_3/\varphi_2 &\\
\cline{1-1}\cline{3-3}\cline{5-5}\cline{7-7}
b_1+r_1&\!+\!&b_2+r_2-r_0\varphi_2/\varphi_1&\!+\!
&b_3+r_3-r_1\varphi_3/\varphi_2&
\!+\cdots.\nonumber
\end{array}
\end{align*}
\end{lemma}
\begin{remark}
This transformation is in particular useful if $b_0+\k\big(a_m/b_m\big)$ and its Bauer-Muir transformation converge to the same value. This was proved to be the case for positive continued fractions and positive modifying factors $r_m$ by Perron~\cite[p.~27]{Perron}, but of course it holds whenever $\k\big(a_m/b_m\big)$ converges with exceptional sequence $\{r_m^\dag\}$ and $\lim\inf \mathfrak{m}(r_m,r_m^\dag)>0$, see \cite[p.~83]{LW1} or \cite{Loren94c}. Berndt~\cite[p.~37-38]{Berndt98} describes this transformation in details.
\end{remark}
\begin{remark}
Although the Bauer-Muir transformation is classical, there is still some room for
further interesting applications. The new idea of this paper is that we always choose suitable nonlinear modifying factors $\left\{r_m\right\}_{m=0}^{\infty}$ so that $\varphi_{m+1}/\varphi_m$ equals to $-1$ for $m\in \mathbb{N}$. In other words, for each proof of five theorems $\varphi_m$ is not constant, which will be called as the \emph{adjoint factors} of $\{r_m\}_{m=0}^{\infty}$ for $b_0+\k\big(a_m/b_m\big)$, or the \emph{adjoint factors} $\left\{\varphi_m\right\}_{m=1}^{\infty}$ for short. In addition, we shall employ the Bauer-Muir transformation twice with different modifying factors. To the best knowledge of authors, such an approach is used for the first time in this paper. So our theorems may be viewed as new applications of the Bauer-Muir transformation. For such approaches, we will show more examples elsewhere.
\end{remark}

\section{The proof of Theorem 1}
For each $j=1,2,3$, we note that the second continued fraction in~\eqref{Theorem 1-Second} is the even part of the first continued fraction in~\eqref{Theorem 1-Second}. So it suffices to prove that the first identity holds in ~\eqref{Theorem 1-First}. In addition, we only guess that the functions in~\eqref{Theorem 1-First} equals the second continued fraction in~\eqref{Theorem 1-Second} by the MC-algorithm.
\subsection{The proof of Theorem 1 for $j=1$}
We shall prove Theorem 1 for $0<\alpha,\beta<1$ and $x>2$. The parabola theorem of Jacobsen~\cite[p.~419--420, Theorem 2.3(iv)]{Jac} can then be employed to extend the domains of convergence for $x, \alpha$ and $\beta$ to those indicated.

Firstly, let us consider the following continued fraction
\begin{equation}
F(x)=\frac{b_0(1;x)}{2}+\K_{m=1}^{\infty}\Bigl(\frac{a_m(1)}{b_m(1;x)}\Bigr),\label{Theorem 1-1-F-D}
\end{equation}
where
\begin{align}
&\begin{cases}
a_{2m}(1)=(2m-2+\alpha+\beta)(2m-1-\alpha)(2m-1-\beta)/4,\\
a_{2m-1}(1)=(2m-\alpha-\beta)(2m-1+\alpha)(2m-1+\beta)/4,
\end{cases}\\
&\begin{cases}
b_{2m}(1;x)=x^2-1,\\
b_{2m-1}(1;x)=2m-1.
\end{cases}
\end{align}
We assume that the adjoint factors $\{\phi_m\}_{m=1}^{\infty}$ are defined by
\begin{equation}\begin{cases}
\phi_{2m}=a_{2m}(1)-r_{2m-1}\left(b_{2m}(1;x)+r_{2m}\right),\\
\phi_{2m-1}=a_{2m-1}(1)-r_{2m-2}\left(b_{2m-1}(1;x)+r_{2m-1}\right),
\end{cases}
\end{equation}
where modifying factors $\{r_m\}_{m=0}^{\infty}$ are of the form of $r_{2m}=u_1 m^2+v_1 m+w_1, r_{2m-1}=u_2 m^2+v_2 m+w_2$.

To prove Theorem 1 for $j=1$ by the Bauer-Muir transformation, we shall choose the parameters $(u_1,u_2,v_1,v_2,w_1,w_2)$ such that both the adjoint subsequence $\phi_{2m}$ and $\phi_{2m-1}$ are constant, i.e. $[m^k]\phi_{2m}=0$ and $[m^k]\phi_{2m-1}=0$ for $k=1,2,3,4$. In a sequel, we shall always use this idea.  With the help of \emph{Mathematica} software, we choose\footnote{We may use another solution $\left(u_1,u_2,v_1,v_2,w_1,w_2\right)=\Bigl(0,-\frac{2 }{-1+x},1-x,\frac{3-x}{-1+x},\frac{1}{2}(1-x^2),\frac{-4x+\rho}{2(1-x)} \Bigr)$ to provide second approach.}
\begin{equation}
\left(u_1,u_2,v_1,v_2,w_1,w_2\right)=\Bigl(0,\frac{2 }{1+x},1+x,-\frac{3+x}{1+x},\frac{1}{2}(1-x^2),\frac{\rho}{2(1+x)} \Bigr),
\end{equation}
where $\rho=2+2x+x^2+\alpha-\alpha^2+\beta-\alpha\beta-\beta^2$. In this case, it is not difficult to check that the following relations hold
\begin{equation}
\phi_{2m}=-\frac 14(x+\alpha)(x+\beta)(x+1-\alpha-\beta)=-\phi_{2m-1}\ne 0.
\end{equation}
We evidently observe that the identity $\phi_{m+1}/\phi_{m}=-1$ always holds for any $m\in\mathbb{N}$. Moreover, $b_0(1;x)/2+r_0=0,\quad b_1(1;x)+r_1=w_2,$ and
\begin{align}\begin{cases}
b_{2m}(1;x)+r_{2m}+r_{2m-2}=(-1+2m)(1+x),\\
 b_{2m+1}(1;x)+r_{2m+1}+r_{2m-1}=\frac{4m^2+\rho}{1+x}.
 \end{cases}
\end{align}
It follows from the Bauer-Muir transformation and an equivalence transformation~(see~\cite[p.~73, Theorem 9]{LW}) that
\begin{align}
F(x)=&\frac{b_0(1;x)}{2}+r_0+\label{Theorem 1-1-RL-1}\\
&\begin{array}{ccccccc}
\phi_1& & -a_1(1)& & -a_2(1) &\\
\cline{1-1}\cline{3-3}\cline{5-5}\cline{7-7}
b_1(1;x)+ r_1&\!+\!&b_2(1;x)+r_2+r_0&\!+\!&b_3(1;x)+r_3+r_1&+\cdots
\end{array}\nonumber\\
=&\begin{array}{ccccccc}
\phi_1& &-a_1(1)& &-a_2(1) &\\
\cline{1-1}\cline{3-3}\cline{5-5}\cline{7-7}
\frac{\rho}{2(1+x)}&\!+\!&1\cdot(1+x)\!&\!+\!&\frac{4+\rho}{1+x}&\!+\!\cdots\!+
\end{array}\nonumber\\
&\begin{array}{ccccc}
-a_{2m-1}(1)&& -a_{2m}(1)&\\
\cline{1-1}\cline{3-3}\cline{5-5}
(2m-1)(1+x)&\!+\!&\frac{4m^2+\rho}{1+x}&\!+\cdots
\end{array}\nonumber\\
=&\begin{array}{ccccccccccc}
-\phi_1(1+x)& & a_1(1)& & a_2(1) &\\
\cline{1-1}\cline{3-3}\cline{5-5}\cline{7-7}
-\rho/2&\!+\!&1&\!+\!&-(4+\rho)&\!+\!\cdots\!+
\end{array}\nonumber\\
&\begin{array}{ccccc}
a_{2m-1}(1)&& a_{2m}(1)&\\
\cline{1-1}\cline{3-3}\cline{5-5}
2m-1&\!+\!&-(4m^2+\rho)&\!+\cdots
\end{array}\nonumber\\
:=&\frac{-\phi_1(1+x)}{-\rho/2+}\K_{m=1}^{\infty}\Bigl(\frac{a_m(1)}{d_m(1;x)}\Bigr),\quad (\mbox{write}).\nonumber
\end{align}
Secondly, to treat the continued fraction $-\rho/2+\k_{m=1}^{\infty}\bigl(a_m(1)/d_m(1;x)\bigr)$, we need to employ the Bauer-Muir transformation again. To do that, we let the adjoint factors have the following form
\begin{equation}\begin{cases}
\psi_{2m}=a_{2m}(1)-\gamma_{2m-1}\left(d_{2m}(1;x)+\gamma_{2m}\right),\\
\psi_{2m-1}=a_{2m-1}(1)-\gamma_{2m-2}\left(d_{2m-1}(1;x)+\gamma_{2m-1}\right),
\end{cases}
\end{equation}
where modifying factors $\gamma_{2m}=p_1 m^2+q_1 m+t_1, \gamma_{2m-1}=p_2 m^2+q_2 m+t_2$.  On taking
\begin{equation}
(p_1,p_2,q_1,q_2,t_1,t_2)=\Bigl(2,0,1-x,-1,\frac {\rho}{2},\frac{3+x}{2}\Bigr),
\end{equation}
we find that
\begin{align*}\begin{cases}
\psi_{2m}=\frac 14(2+x-\alpha)(2+x-\beta)(1+x+\alpha+\beta)=-\psi_{2m-1}\ne 0,\\
d_{2m}(1;x)+\gamma_{2m}+\gamma_{2m-2}=-(2m-1)(1+x),\\
 d_{2m+1}(1;x)+\gamma_{2m+1}+\gamma_{2m-1}=3+x.
\end{cases}\end{align*}
We see easily that the identity $\psi_{m+1}/\psi_{m}=-1$ always holds for any $m\in\mathbb{N}$. By using the  Bauer-Muir transformation again, and then employing an equivalence transformation, one derives that
\begin{align}
&-\frac {\rho}{2}+\K_{m=1}^{\infty}\Bigl(\frac{a_m(1)}{d_m(1;x)}\Bigr)\label{Theorem 1-1-RL-2}\\
=&-\frac {\rho}{2}+\gamma_0+\nonumber\\
&\begin{array}{cccccccc}
\psi_1& & -a_1(1)& & -a_2(1) &\\
\cline{1-1}\cline{3-3}\cline{5-5}\cline{7-7}
d_1(1;x)+ \gamma_1&\!+\!&d_2(1;x)+\gamma_2+\gamma_0&\!+\!&d_3(1;x)+\gamma_3+\gamma_1&\!+\cdots
\end{array}\nonumber\\
=&\begin{array}{ccccccc}
\psi_1& &-a_1(1)& &-a_2(1) &\\
\cline{1-1}\cline{3-3}\cline{5-5}
(3+x)/2&\!+\!&-1\cdot(1+x)&\!+\!&3+x&\!+\!\cdots\!+
\end{array}\nonumber\\
&\begin{array}{ccccc}
-a_{2m-1}(1)&& -a_{2m}(1)&\\
\cline{1-1}\cline{3-3}\cline{5-5}
-(2m-1)(1+x)&\!+\!&3+x&\!+\cdots
\end{array}\nonumber\\
=&\begin{array}{ccccccccccc}
\psi_1& &a_1(1)& &a_2(1) &\\
\cline{1-1}\cline{3-3}\cline{5-5}
(3+x)/2&\!+\!&1\cdot(1+x)&\!+\!&3+x&\!+\!\cdots\!+
\end{array}\nonumber\\
&\begin{array}{ccccc}
a_{2m-1}(1)&& a_{2m}(1)&\\
\cline{1-1}\cline{3-3}\cline{5-5}
(2m-1)(1+x)&\!+\!&3+x&\!+\cdots
\end{array}\nonumber\\
=&\begin{array}{ccccccccccc}
\psi_1(1+x)& &a_1(1)& &a_2(1) &\\
\cline{1-1}\cline{3-3}\cline{5-5}
(3+x)(1+x)/2&\!+\!&1&\!+\!&(3+x)(1+x)&\!+\!\cdots\!+
\end{array}\nonumber\\
&\begin{array}{ccccc}
a_{2m-1}(1)&& a_{2m}(1)&\\
\cline{1-1}\cline{3-3}\cline{5-5}
2m-1&\!+\!&(3+x)(1+x)&\!+\cdots\!
\end{array}\nonumber\\
=&\frac{\psi_1(1+x)}{F(x+2)}.\nonumber
\end{align}
Thirdly, substituting \eqref{Theorem 1-1-RL-2} into \eqref{Theorem 1-1-RL-1},  after some simplifications we deduce that
\begin{equation}
\frac{F(x)}{F(x+2)}=\frac{\left(x+\alpha\right)\left(x+1-\alpha-\beta\right)\left(x+\beta\right)}
{\left(x+2-\alpha\right)(x+1+\alpha+\beta)\left(x+2-\beta\right)}.\label{Theorem 1-Recurrence Relation}
\end{equation}
By iteration of the foregoing identity, we find that, for each $m\in\mathbb{N}$
\begin{align}
&\frac{F(x)}{F(x+2m)}\\
=&\frac{\prod_{k=0}^{m-1}\left(x+2k+\alpha\right)\left(x+2k+1-\alpha-\beta\right)\left(x+2k+\beta\right)}
{\prod_{k=0}^{m-1}\left(x+2k+2-\alpha\right)(x+2k+1+\alpha+\beta)\left(x+2k+2-\beta\right)}\nonumber
\\
=&\frac{\prod_{k=0}^{m-1}\left(\frac 12(x+\alpha)+k\right)\left(\frac 12(x+1-\alpha-\beta)+k\right)\left(\frac 12(x+\beta)+k\right)}{\prod_{k=0}^{m-1}\left(\frac 12(x+2-\alpha)+k\right)\left(\frac 12(x+1+\alpha+\beta)+k\right)\left(\frac 12(x+2-\beta)+k\right)}\nonumber\\
=&\frac {1}{m^2} \frac{\prod_{k=0}^{m-1}\left(\frac 12(x+\alpha)+k\right)(m!)^{-1}m^{-\frac 12(x+\alpha)}}{\prod_{k=0}^{m-1}\left(\frac 12(x+2-\alpha)+k\right)(m!)^{-1}m^{-\frac 12(x+2-\alpha)}}
\nonumber\\
&\times \frac{\prod_{k=0}^{m-1}\left(\frac 12(x+1-\alpha-\beta)+k\right)(m!)^{-1}m^{-\frac 12(x+1-\alpha-\beta)}}{\prod_{k=0}^{m-1}\left(\frac 12(x+1+\alpha+\beta)+k\right)(m!)^{-1}m^{-\frac 12(x+1+\alpha+\beta)}}\nonumber\\
&\times\frac{\prod_{k=0}^{m-1}\left(\frac 12(x+\beta)+k\right)(m!)^{-1}m^{-\frac 12(x+\beta)}}{\prod_{k=0}^{m-1}\left(\frac 12(x+2-\beta)+k\right)(m!)^{-1}m^{-\frac 12(x+2-\beta)}}.\nonumber
\end{align}
By means of Euler's formula~(see~\cite[p.~255, (6.1.2)]{AS})
\begin{equation}\label{Euler's formula}
\Gamma(z)=\lim_{n\rightarrow\infty}\frac{n!~n^z}{z(z+1)\cdots(z+n)}\quad(z\neq 0,-1,-2,\ldots),
\end{equation}
we derive that
\begin{align}
&\lim_{m\rightarrow\infty}\frac{F(x)m^2}{F(x+2m)}\label{Theorem 1-1-3}\\
=&\frac{
\Gamma\left(\frac 12 ( x+2 - \alpha)\right) \Gamma\left(\frac 12( x+1+\alpha+\beta)\right)\Gamma
 \left( \frac 12 (x +2-\beta)\right)}{ \Gamma\left(\frac 12 ( x+\alpha)\right) \Gamma\left(\frac 12(x+1-\alpha-\beta)\right)\Gamma
  \left(\frac 12 (x +\beta)\right)}.\nonumber
\end{align}
From the definition of $F(x)$ in \eqref{Theorem 1-1-F-D}, it is not difficult to prove that
\begin{equation}
\lim_{m\rightarrow\infty}\frac{m^2}{F(x+2m)}=\frac 12.\label{Theorem 1-1-4}
\end{equation}
Finally, the first identity in \eqref{Theorem 1-First} follows from \eqref{Theorem 1-1-F-D}, \eqref{Theorem 1-1-3}, \eqref{Theorem 1-1-4} and an equivalence transformation readily, and this completes the proof of Theorem 1 for $j=1$.\qed.

\subsection{The proof of Theorem 1 for $j=2$}
Firstly, let us consider the following continued fraction
\begin{equation}
F(x)=\frac{b_0(2;x)}{2}+\K_{m=1}^{\infty}\Bigl(\frac{a_m(2)}{b_m(2;x)}\Bigr),\label{Theorem 1-2-F-D}
\end{equation}
where
\begin{align}
&\begin{cases}
a_{2m}(2)=(2m-2+\alpha+\beta)(2m-1+\alpha)(2m-1+\beta)/4,\\
a_{2m-1}(2)=(2m-\alpha-\beta)(2m-1-\alpha)(2m-1-\beta)/4,
\end{cases}\\
&\begin{cases}
b_{2m}(2;x)=x^2-(1-\alpha-\beta)^2,\\
b_{2m-1}(2;x)=2m-1.
\end{cases}
\end{align}
Throughout this subsection, for brevity, we use $\omega$ to denote $x^2-(1-\alpha-\beta)^2$. We let
\begin{equation}\begin{cases}
\phi_{2m}=a_{2m}(2)-r_{2m-1}\left(b_{2m}(2;x)+r_{2m}\right),\\
\phi_{2m-1}=a_{2m-1}(2)-r_{2m-2}\left(b_{2m-1}(2;x)+r_{2m-1}\right),
\end{cases}
\end{equation}
where $r_{2m}=u_1 m^2+v_1 m+w_1, r_{2m-1}=u_2 m^2+v_2 m+w_2$.

To prove Theorem 1 for $j=2$, we choose
\begin{equation}
\left(u_1,u_2,v_1,v_2,w_1,w_2\right)=\Bigl(0,\frac{2}{\rho}, \rho,-1-\frac{2}{\rho},-\frac{\omega}{2},\frac{h}{2\rho} \Bigr),
\end{equation}
where $\rho=x+1-\alpha-\beta$ and $h=(x+1)^2+1-\alpha-\beta+\alpha\beta$. In this case, it is not difficult to check that the following relations hold
\begin{equation}
\phi_{2m}=-\frac 14(x+1)(x+\alpha)(x+\beta)=-\phi_{2m-1}\ne 0.
\end{equation}
Moreover, $b_0(2;x)/2+r_0=0,\quad b_1(2;x)+r_1=w_2,$ and
\begin{align}\begin{cases}
b_{2m}(2;x)+r_{2m}+r_{2m-2}=(2m-1)\rho,\\ b_{2m+1}(2;x)+r_{2m+1}+r_{2m-1}=\frac{4m^2+h}{\rho}.
\end{cases}
\end{align}
It follows from the Bauer-Muir transformation and an equivalence transformation that
\begin{align}
F(x)=&\frac{b_0(2;x)}{2}+r_0+\label{Theorem 1-2-RL-1}\\
&\begin{array}{ccccccc}
\phi_1& & -a_1(2)& & -a_2(2) &\\
\cline{1-1}\cline{3-3}\cline{5-5}\cline{7-7}
b_1(2;x)+r_1&\!+\!&b_2(2;x)+r_2+r_0&\!+\!&b_3(2;x)+r_3+r_1&\!+\!\cdots
\end{array}\nonumber\\
=&\begin{array}{ccccccc}
\phi_1& &-a_1(2)& &-a_2(2) &\\
\cline{1-1}\cline{3-3}\cline{5-5}
w_2&\!+\!&1\cdot\rho&\!+\!&(4+h)/\rho&\!+\!\cdots\!+
\end{array}\nonumber\\
&\begin{array}{ccccc}
-a_{2m-1}(2)&& -a_{2m}(2)&\\
\cline{1-1}\cline{3-3}\cline{5-5}
(2m-1)\rho&\!+\!&(4m^2+h)/\rho&\!+\cdots
\end{array}\nonumber\\
=&\begin{array}{ccccccc}
-\phi_1\rho& & a_1(2)& & a_2(2) &\\
\cline{1-1}\cline{3-3}\cline{5-5}
-h/2&\!+\!&1&\!+\!&-(4+h)&\!+\!\cdots\!+
\end{array}\nonumber\\
&\begin{array}{ccccc}
a_{2m-1}(2)&& a_{2m}(2)&\\
\cline{1-1}\cline{3-3}\cline{5-5}
2m-1&\!+\!&-(4m^2+h)&\!+\cdots
\end{array}\nonumber\\
:=&\frac{\phi_1\rho}{-h/2+}\K_{m=1}^{\infty}\Bigl(\frac{a_m(2)}{d_m(2;x)}\Bigr),\quad (\mbox{write}).\nonumber
\end{align}
Secondly, to treat the continued fraction $-h/2+\k_{m=1}^{\infty}\bigl(a_m(2)/d_m(2;x)\bigr)$, we need to employ the Bauer-Muir transformation again. To this end, we let
\begin{equation}\begin{cases}
\psi_{2m}=a_{2m}(2)-\gamma_{2m-1}\left(d_{2m}(2;x)+\gamma_{2m}\right),\\
\psi_{2m-1}=a_{2m-1}(2)-\gamma_{2m-2}\left(d_{2m-1}(2;x)+\gamma_{2m-1}\right),
\end{cases}
\end{equation}
where $\gamma_{2m}=p_1 m^2+q_1 m+t_1, \gamma_{2m-1}=p_2 m^2+q_2 m+t_2$. We write $\sigma=x-1+\alpha+\beta$. On taking
\begin{equation}
(p_1,p_2,q_1,q_2,t_1,t_2)=\Bigl(2,0,-\sigma,-1,\frac h2,1+\frac{\rho}{2}\Bigr),
\end{equation}
we find that
\begin{align}\begin{cases}
\psi_{2m}=\frac 14(x+1)(x+2-\alpha)(x+2-\beta)=-\psi_{2m-1}\ne 0,\\
d_{2m}(2;x)+\gamma_{2m}+\gamma_{2m-2}=-(2m-1)(\sigma+2),\\
d_{2m+1}(2;x)+\gamma_{2m+1}+\gamma_{2m-1}=\rho+2.
\end{cases}\end{align}
By using the  Bauer-Muir transformation again, and then employing an equivalence transformation, one derives that
\begin{align}
&-\frac h2+\K_{m=1}^{\infty}\Bigl(\frac{a_m(2)}{d_m(2;x)}\Bigr)\label{Theorem 1-2-RL-2}\\
=&-\frac h2+\gamma_0+\nonumber\\
&\begin{array}{cccccccc}
\psi_1& & -a_1(2)& & -a_2(2) &\\
\cline{1-1}\cline{3-3}\cline{5-5}\cline{7-7}
d_1(2;x)+\gamma_1&\!+\!&d_2(2;x)+\gamma_2+\gamma_0&\!+\!&d_3(2;x)+\gamma_3+\gamma_1&\!+\!\cdots
\end{array}\nonumber\\
=&\begin{array}{ccccccc}
\psi_1& &-a_1(2)& &-a_2(2) &\\
\cline{1-1}\cline{3-3}\cline{5-5}
(\rho+2)/2&\!+\!&-1\cdot(\sigma+2)&\!+\!&\rho+2&\!+\!\cdots\!+
\end{array}\nonumber\\
&\begin{array}{ccccc}
-a_{2m-1}(2)&& -a_{2m}(2)&\\
\cline{1-1}\cline{3-3}\cline{5-5}
-(2m-1)(\sigma+2)&\!+\!&\rho+2&\!+\cdots\!
\end{array}\nonumber\\
=&\begin{array}{ccccccc}
\psi_1& &a_1(2)& &a_2(2) &\\
\cline{1-1}\cline{3-3}\cline{5-5}
(\rho+2)/2&\!+\!&1\cdot(\sigma+2)&\!+\!&\rho+2&\!+\!\cdots\!+
\end{array}\nonumber\\
&\begin{array}{ccccc}
a_{2m-1}(2)&& a_{2m}(2)&\\
\cline{1-1}\cline{3-3}\cline{5-5}
(2m-1)(\sigma+2)&\!+\!&\rho+2&\!+\cdots
\end{array}\nonumber\\
=&\begin{array}{ccccccc}
\psi_1(\sigma+2)& &a_1(2)& &a_2(2) &\\
\cline{1-1}\cline{3-3}\cline{5-5}
(\rho+2)(\sigma+2)/2&\!+\!&1&\!+\!&(\rho+2)(\sigma+2)&\!+\!\cdots\!+
\end{array}
\nonumber\\
&\begin{array}{ccccc}
a_{2m-1}(2)&& a_{2m}(2)&\\
\cline{1-1}\cline{3-3}\cline{5-5}
2m-1&\!+\!&(\rho+2)(\sigma+2)&\!+\cdots
\end{array}\nonumber\\
=&\frac{\psi_1(\sigma+2)}{F(x+2)}.\nonumber
\end{align}
Thirdly, substituting \eqref{Theorem 1-2-RL-2} into \eqref{Theorem 1-2-RL-1}, after some simplifications we deduce that
\eqref{Theorem 1-Recurrence Relation} still holds. Just as the proof of the previous case, one may completes the proof of Theorem 1 for $j=2$.\qed

\subsection{The proof of Theorem 1 for $j=3$}

Firstly, let us consider the following continued fraction
\begin{equation}
F(x)=\frac{b_0(3;x)}{2}+\K_{m=1}^{\infty}\Bigl(\frac{a_m(3)}{b_m(3;x)}\Bigr),\label{Theorem 1-3-F-D}
\end{equation}
where
\begin{align}
&\begin{cases}
a_{2m}(3)=(2m-\alpha-\beta)(2m-1+\alpha)(2m-1-\beta)/4,\\
a_{2m-1}(3)=(2m-2+\alpha+\beta)(2m-1-\alpha)(2m-1+\beta)/4,
\end{cases}\\
&\begin{cases}
b_{2m}(3;x)=x^2-\beta^2,\\
b_{2m-1}(3;x)=2m-1.
\end{cases}
\end{align}
We let
\begin{equation}\begin{cases}
\phi_{2m}=a_{2m}(3)-r_{2m-1}\left(b_{2m}(3;x)+r_{2m}\right),\\
\phi_{2m-1}=a_{2m-1}(3)-r_{2m-2}\left(b_{2m-1}(3;x)+r_{2m-1}\right),
\end{cases}
\end{equation}
where $r_{2m}=u_1 m^2+v_1 m+w_1, r_{2m-1}=u_2 m^2+v_2 m+w_2$.

To prove Theorem 1 for $j=3$, we choose
\begin{align}
&\left(u_1,u_2,v_1,v_2,w_1,w_2\right)\\
=&\Bigl(0,\frac{2 }{x+\beta},x+\beta,-\frac{2+x+\beta}{x+\beta},\frac{1}{2}(-x^2+\beta^2),\frac{\rho}{2(x+\beta)} \Bigr),\nonumber
\end{align}
where $\rho=1+2x+x^2+\alpha-\alpha^2+\beta-\alpha\beta$. In this case, it is not difficult to check that the following relations hold
\begin{equation}
\phi_{2m}=-\frac 14(x+1)(x+\alpha)(x+1-\alpha-\beta)=-\phi_{2m-1}\ne 0.
\end{equation}
Moreover, $b_0(3;x)/2+r_0=0,\quad b_1(3;x)+r_1=w_2,$ and
\begin{align}\begin{cases}
b_{2m}(3;x)+r_{2m}+r_{2m-2}=(-1+2m)(x+\beta),\\
 b_{2m+1}(3;x)+r_{2m+1}+r_{2m-1}=\frac{4m^2+\rho}{x+\beta}.
 \end{cases}
\end{align}
It follows from the Bauer-Muir transformation and an equivalence transformation that
\begin{align}
F(x)=&\frac{b_0(3;x)}{2}+r_0+\label{Theorem 1-3-RL-1}\\
&\begin{array}{cccccccc}
\phi_1& & -a_1(3)& & -a_2(3) &\\
\cline{1-1}\cline{3-3}\cline{5-5}\cline{7-7}
b_1+ r_1&\!+\!&b_2+r_2+r_0&\!+\!&b_3+r_3+r_1&+\cdots
\end{array}\nonumber\\
=&\begin{array}{ccccccccccc}
\phi_1& &-a_1(3)& &-a_2(3) &\\
\cline{1-1}\cline{3-3}\cline{5-5}
w_2&\!+\!&1\cdot(x+\beta)&\!+\!&(4\cdot 1^2+\rho)/(x+\beta)&\!+\!\cdots\!+
\end{array}\nonumber\\
&\begin{array}{ccccc}
-a_{2m-1}(3)&& -a_{2m}(3)&\\
\cline{1-1}\cline{3-3}\cline{5-5}
(2m-1)(x+\beta)&\!+\!&(4m^2+\rho)/(x+\beta)&\!+\cdots
\end{array}\nonumber\\
=&\begin{array}{ccccccc}
-\phi_1(x+\!\beta)& & a_1(3)& & a_2(3) &\\
\cline{1-1}\cline{3-3}\cline{5-5}
-\rho/2&\!+\!&1&\!+\!&-(4+\rho)&\!+\!\cdots\!+
\end{array}\nonumber\\
&\begin{array}{ccccc}
a_{2m-1}(3)&& a_{2m}(3)&\\
\cline{1-1}\cline{3-3}\cline{5-5}
2m-1&\!+\!&-(4m^2+\rho)&\!+\cdots
\end{array}\nonumber\\
:=&\frac{-\phi_1(x+\beta)}{-\rho/2+}\K_{m=1}^{\infty}\Bigl(\frac{a_m(3)}{d_m(3;x)}\Bigr),\quad (\mbox{write}).\nonumber
\end{align}
Secondly, to treat the continued fraction $-\rho/2+\k_{m=1}^{\infty}\bigl(a_m(3)/d_m(3;x)\bigr)$, we need to employ the Bauer-Muir transformation again. To do that, we let
\begin{equation}\begin{cases}
\psi_{2m}=a_{2m}(3)-\gamma_{2m-1}\left(d_{2m}(3;x)+\gamma_{2m}\right),\\
\psi_{2m-1}=a_{2m-1}(3)-\gamma_{2m-2}\left(d_{2m-1}(3;x)+\gamma_{2m-1}\right),
\end{cases}
\end{equation}
where $\gamma_{2m}=p_1 m^2+q_1 m+t_1, \gamma_{2m-1}=p_2 m^2+q_2 m+t_2$.  On taking
\begin{equation}
(p_1,p_2,q_1,q_2,t_1,t_2)=\Bigl(2,0,-x+\beta,-1,\frac {\rho}{2},\frac{2+x+\beta}{2}\Bigr),
\end{equation}
we find that
\begin{align}\begin{cases}
\psi_{2m}=\frac 14(1+x)(2+x-\alpha)(1+x+\alpha+\beta)=-\psi_{2m-1}\ne 0,\\
d_{2m}(3;x)+\gamma_{2m}+\gamma_{2m-2}=-(2m-1)(2+x-\beta),\\
d_{2m+1}(3;x)+\gamma_{2m+1}+\gamma_{2m-1}=2+x+\beta.
\end{cases}\end{align}
By using the Bauer-Muir transformation again, and then employing an equivalence transformation, one derives that
\begin{align}
&-\frac {\rho}{2}+\K_{m=1}^{\infty}\Bigl(\frac{a_m(3)}{d_m(3;x)}\Bigr)\nonumber\\
=&-\frac {\rho}{2}+\gamma_0+\nonumber\\
&\begin{array}{ccccccccc}
\psi_1& & -a_1(3)& & -a_2(3) &\\
\cline{1-1}\cline{3-3}\cline{5-5}\cline{7-7}
d_1+ \gamma_1&\!+\!&d_2+\gamma_2+\gamma_0&\!+\!&d_3+\gamma_3+\gamma_1&\!+\cdots
\end{array}\nonumber\\
=&\begin{array}{ccccccc}
\psi_1& &-a_1(3)& &-a_2(3) &\\
\cline{1-1}\cline{3-3}\cline{5-5}
(2+x+\beta)/2&\!+\!&-1\cdot(2+x-\beta)&\!+\!&2+x+\beta&\!+\!\cdots\!+
\end{array}\nonumber\\
&\begin{array}{ccccc}
-a_{2m-1}(3)&&-a_{2m}(3)&\\
\cline{1-1}\cline{3-3}\cline{5-5}
-(2m-1)(2+x-\beta)&\!+\!&2+x+\beta&\!+\!\cdots\!
\end{array}\nonumber\\
=&\begin{array}{ccccccc}
\psi_1& &a_1(3)& &a_2(3) &\\
\cline{1-1}\cline{3-3}\cline{5-5}
(2+x+\beta)/2&\!+\!&1\cdot(2+x-\beta)&\!+\!&2+x+\beta&\!+\!\cdots\!+
\end{array}\nonumber\\
&\begin{array}{ccccc}
a_{2m-1}(3)&& a_{2m}(3)&\\
\cline{1-1}\cline{3-3}\cline{5-5}
(2m-1)(2+x-\beta)&\!+\!&2+x+\beta&\!+\!\cdots
\end{array}\nonumber\\
=&\begin{array}{ccccccc}
\psi_1(2+x-\beta)& &a_1(3)& &a_2(3) &\\
\cline{1-1}\cline{3-3}\cline{5-5}
(2+x+\beta)(2+x-\beta)/2&\!+\!&1&\!+\!&\sigma&\!+\!\cdots\!+
\end{array}\nonumber\\
&\begin{array}{ccccc}
a_{2m-1}(3)&& a_{2m}(3)&\\
\cline{1-1}\cline{3-3}\cline{5-5}
2m-1&\!+\!&\sigma&\!+\cdots\!
\end{array}\nonumber\\
=&\frac{\psi_1(2+x-\beta)}{F(x+2)},\nonumber
\end{align}
where $\sigma=(2+x+\beta)(2+x-\beta)$. Thirdly, substituting the foregoing formula 
into \eqref{Theorem 1-3-RL-1},  after some simplifications we deduce that \eqref{Theorem 1-Recurrence Relation} also holds. Similarly for case $j=1$, we shall finish the proof of Theorem 1 for $j=3$.\qed

\section{The proof of Theorem 2}
\proof  We shall prove Theorem 2 for $0<\alpha,\beta<1$ and $x>2$. The parabola theorem of Jacobsen~\cite[p.~419--420, Theorem 2.3(iv)]{Jac} can then be employed to extend the domains of convergence for $x, \alpha$ and $\beta$ to those indicated.

Firstly, let us consider the following continued fraction
\begin{equation}
F(x)=\frac{b_0(x)}{2}+\K_{m=1}^{\infty}\Bigl(\frac{a_m}{b_m(x)}\Bigr),\label{Theorem 2-F-D}
\end{equation}
where
\begin{align}
&\begin{cases}
a_{2m}=(2m-\alpha-\beta)(2m-1+\alpha)(2m-1+\beta)/4,\\
a_{2m-1}=(2m-2+\alpha+\beta)(2m-1-\alpha)(2m-1-\beta)/4,
\end{cases}\\
&\begin{cases}
b_{2m}(x)=x,\\
b_{2m-1}(x)=(2m-1)x.
\end{cases}
\end{align}
Suppose that the adjoint factors $\{\phi_m\}_{m=1}^{\infty}$ are defined by
\begin{equation}\begin{cases}
\phi_{2m}=a_{2m}-r_{2m-1}\left(b_{2m}(x)+r_{2m}\right),\\
\phi_{2m-1}=a_{2m-1}-r_{2m-2}\left(b_{2m-1}(x)+r_{2m-1}\right),
\end{cases}
\end{equation}
where modifying factors $\{r_m\}_{m=0}^{\infty}$ have the form of $r_{2m}=u_1 m^2+v_1 m+w_1, r_{2m-1}=u_2 m^2+v_2 m+w_2$.

To prove Theorem 2 by the Bauer-Muir transformation, we shall choose the parameters $(u_1,u_2,v_1,v_2,w_1,w_2)$ such that both $\phi_{2m}$ and $\phi_{2m-1}$ are constant, i.e. $[m^j]\phi_{2m}=0$ and $[m^j]\phi_{2m-1}=0$ for $j=1,2,3,4$. With the aid of \emph{Mathematica} software, let us choose
\begin{align}
&\left(u_1,u_2,v_1,v_2,w_1,w_2\right)\\
=&\Bigl(0,-2,-1,2-x,-\frac x2,\frac 12(-(x-1)^2+\alpha^2-\alpha+\alpha\beta+\beta^2-\beta)\Bigr).\nonumber
\end{align}
In this case, it is not difficult to check that the following relations hold
\begin{equation}
\phi_{2m}=\frac 14(x-1+\alpha)(x-\alpha-\beta)(x-1+\beta)=-\phi_{2m-1}\ne 0.
\end{equation}
We evidently observe that the identity $\phi_{m+1}/\phi_{m}=-1$ always holds for any $m\in\mathbb{N}$. Moreover, $b_0/2+r_0=0,\quad b_1+r_1=w_2,$ and
\begin{align}\begin{cases}
b_{2m}(x)+r_{2m}+r_{2m-2}=-2m+1,\\ b_{2m+1}(x)+r_{2m+1}+r_{2m-1}=-4m^2+2w_2.
\end{cases}
\end{align}
It follows from the Bauer-Muir transformation and an equivalence transformation~(see~\cite[p.~73, Theorem 9]{LW}) that
\begin{align}
F(x)=&\frac{b_0(x)}{2}+r_0+\label{Theorem 2-RL-1}\\
&\begin{array}{ccccccc}
\phi_1& & -a_1& & -a_2 &\\
\cline{1-1}\cline{3-3}\cline{5-5}\cline{7-7}
b_1(x)+ r_1&\!+\!&b_2(x)+r_2+r_0&\!+\!&b_3(x)+r_3+r_1&\!+\cdots
\end{array}\nonumber\\
=&\begin{array}{ccccccccccc}
\phi_1& &-a_1& &-a_2 &&\\
\cline{1-1}\cline{3-3}\cline{5-5}
w_2&\!+\!&-1&\!+\!&-4+2w_2&\!+\!\cdots\!+
\end{array}\nonumber\\
&\begin{array}{cccccc}
-a_{2m-1}&& -a_{2m}&\\
\cline{1-1}\cline{3-3}\cline{5-5}
-2m+1&\!+\!&-4m^2+2w_2&\!+\cdots
\end{array}\nonumber\\
=&\begin{array}{ccccccc}
\phi_1& & a_1& & a_2 &\\
\cline{1-1}\cline{3-3}\cline{5-5}
w_2&\!+\!&1&\!+\!&-4+2w_2&\!+\!\cdots\!+
\end{array}\nonumber\\
&\begin{array}{cccccc}
a_{2m-1}&& a_{2m}&\\
\cline{1-1}\cline{3-3}\cline{5-5}
2m-1&\!+\!&-4m^2+2w_2&\!+\cdots
\end{array}\nonumber\\
:=&\frac{\phi_1}{w_2+}\K_{m=1}^{\infty}\Bigl(\frac{a_m}{d_m(x)}\Bigr),\quad (\mbox{say}).\nonumber
\end{align}
Secondly, to treat the continued fraction $w_2+\k_{m=1}^{\infty}\bigl(a_m/d_m(x)\bigr)$, we require to employ the Bauer-Muir transformation again. To this end, we let
\begin{equation}\begin{cases}
\psi_{2m}=a_{2m}-\gamma_{2m-1}\left(d_{2m}(x)+\gamma_{2m}\right),\\
\psi_{2m-1}=a_{2m-1}-\gamma_{2m-2}\left(d_{2m-1}(x)+\gamma_{2m-1}\right),
\end{cases}\end{equation}
where $\gamma_{2m}=p_1 m^2+q_1 m+t_1, \gamma_{2m-1}=p_2 m^2+q_2 m+t_2$. Very similar to previous approach, on taking
\begin{equation}
(p_1,p_2,q_1,q_2,t_1,t_2)=\Bigl(2,0,x,-1,-w_2,1-\frac x2\Bigr),
\end{equation}
we find that
\begin{align}\begin{cases}
\psi_{2m}=-\frac 14(x-1-\alpha)(x-2+\alpha+\beta)(x-1-\beta)=-\psi_{2m-1}\ne 0,\\
d_{2m}(x)+\gamma_{2m}+\gamma_{2m-2}=(2m-1)(x-2),\\
d_{2m+1}(x)+\gamma_{2m+1}+\gamma_{2m-1}=-(x-2).
\end{cases}\end{align}
We see easily that the identity $\psi_{m+1}/\psi_{m}=-1$ always holds for any $m\in\mathbb{N}$. Applying the Bauer-Muir transformation and an equivalence transformation again, one derives that
\begin{align}
&w_2+\K_{m=1}^{\infty}\Bigl(\frac{a_m}{d_m(x)}\Bigr)\label{Theorem 2-RL-2}\\
=&w_2+\gamma_0+\begin{array}{cccccccc}
\psi_1& & -a_1& & -a_2 &\\
\cline{1-1}\cline{3-3}\cline{5-5}\cline{7-7}
d_1(x)+ \gamma_1&\!+\!&d_2(x)+\gamma_2+\gamma_0&\!+\!&d_3(x)+\gamma_3+\gamma_1&\!+\cdots
\end{array}\nonumber\\
=&\begin{array}{ccccccc}
\psi_1& &-a_1& &-a_2 &\\
\cline{1-1}\cline{3-3}\cline{5-5}
-(x-2)/2&\!+\!&1\cdot(x-2)&\!+\!&-(x-2)&\!+\!\cdots\!+
\end{array}\nonumber\\
&\begin{array}{cccccc}
-a_{2m-1}&& -a_{2m}&\\
\cline{1-1}\cline{3-3}\cline{5-5}
(2m-1)(x-2)&\!+\!&-(x-2)&\!+\cdots\!
\end{array}\nonumber\\
=&\begin{array}{ccccccc}
-\psi_1& &a_1& &a_2 &\\
\cline{1-1}\cline{3-3}\cline{5-5}
(x-2)/2&\!+\!&1\cdot(x-2)&\!+\!&x-2&\!+\!\cdots\!+
\end{array}\nonumber\\
&\begin{array}{cccccc}
a_{2m-1}&& a_{2m}&\\
\cline{1-1}\cline{3-3}\cline{5-5}
(2m-1)(x-2)&\!+\!&x-2&\!+\cdots
\end{array}\nonumber\\
=&\frac{-\psi_1}{F(x-2)}.\nonumber
\end{align}
Thirdly, substituting \eqref{Theorem 2-RL-2} into \eqref{Theorem 2-RL-1}, we deduce that
\begin{equation}
\frac{F(x)}{F(x-2)}=\frac{-\frac 14(x-1+\alpha)(x-\alpha-\beta)(x-1+\beta)}{-\frac 14(x-1-\alpha)(x-2+\alpha+\beta)(x-1-\beta)}.
\end{equation}
Replacing $x$ by $x+2$ in the above equality and then taking reciprocals, we get
\begin{equation}
\frac{F(x)}{F(x+2)}=\frac{(x+1-\alpha)(x+\alpha+\beta)(x+1-\beta)}{(x+1+\alpha)(x+2-\alpha-\beta)(x+1+\beta)}.
\end{equation}
By iteration of this formula, we find that, for each $m\in\mathbb{N}$
\begin{align}
&\frac{F(x)}{F(x+2m)}\label{Theorem 2-RL-3}\\
=&\frac{\prod_{k=0}^{m-1}(x+2k+1-\alpha)(x+2k+\alpha+\beta)(x+2k+1-\beta)}
{\prod_{k=0}^{m-1}(x+2k+1+\alpha)(x+2k+2-\alpha-\beta)(x+2k+1+\beta)}\nonumber\\
=&\frac{\prod_{k=0}^{m-1}\left(\frac 12(x+1-\alpha)+k\right)\left(\frac 12(x+\alpha+\beta)+k\right)\left(\frac 12(x+1-\beta)+k\right)}
{\prod_{k=0}^{m-1}\left(\frac 12(x+1+\alpha)+k\right)\left(\frac 12(x+2-\alpha-\beta)+k\right)\left(\frac 12(x+1+\beta)+k\right)}
\nonumber\\
=&\frac 1m \frac{\prod_{k=0}^{m-1}\left(\frac 12(x+1-\alpha)+k\right)(m!)^{-1}m^{-\frac 12(x+1-\alpha)}}{\prod_{k=0}^{m-1}\left(\frac 12(x+1+\alpha)+k\right)(m!)^{-1}m^{-\frac 12(x+1+\alpha)}}
\times
\nonumber\\
& \frac{\prod_{k=0}^{m-1}\left(\frac 12(x+\alpha+\beta)+k\right)(m!)^{-1}m^{-\frac 12(x+\alpha+\beta)}}{\prod_{k=0}^{m-1}\left(\frac 12(x\!+\!2\!-\!\alpha\!-\!\beta)\!+\!k\right)(m!)^{-1}m^{-\frac 12(x\!+\!2\!-\!\alpha\!-\!\beta)}}\times\nonumber\\
&\frac{\prod_{k=0}^{m-1}\left(\frac 12(x+1-\beta)+k\right)(m!)^{-1}m^{-\frac 12(x+1-\beta)}}{\prod_{k=0}^{m-1}\left(\frac 12(x+1+\beta)+k\right)(m!)^{-1}m^{-\frac 12(x+1+\beta)}}.\nonumber
\end{align}
By means of Euler's formula \eqref{Euler's formula}, it follows from the identity \eqref{Theorem 2-RL-3} that
\begin{align}
&\lim_{m\rightarrow\infty}\frac{F(x)m}{F(x+2m)}\label{Theorem 2-3}\\
=&\frac{\Gamma\left(\frac 12(x+1+\alpha)\right)\Gamma\left(\frac 12(x+2-\alpha-\beta)\right)\Gamma\left(\frac 12(x+1+\beta)\right)}{\Gamma\left(\frac 12(x+1-\alpha)\right)\Gamma\left(\frac 12(x+\alpha+\beta)\right)\Gamma\left(\frac 12(x+1-\beta)\right)}.\nonumber
\end{align}
From the definition of $F(x)$ in \eqref{Theorem 2-F-D}, we easily show that
\begin{equation}
\lim_{m\rightarrow\infty}\frac{m}{F(x+2m)}=1.\label{Theorem 2-4}
\end{equation}
Finally, combining \eqref{Theorem 2-F-D}, \eqref{Theorem 2-3} and \eqref{Theorem 2-4}, and using an equivalence transformation again, we deduce \eqref{Theorem 2} as a desired result.\qed

\section{The proof of Theorem 3}
\proof  We shall prove Theorem 3 for $0<\alpha,\beta<1$ and $x>2$. An argument of Jacobsen~\cite[p.~419--420, Theorem 2.3(iv)]{Jac} can then be employed to extend the domains of convergence for $x, l, n,$ and $\eta$ to those indicated.

Firstly, let us consider the following continued fraction
\begin{equation}
F(x)=\frac{b_0(x)}{2}+\K_{m=1}^{\infty}\Bigl(\frac{a_m}{b_m(x)}\Bigr),\label{Theorem 3-F-D}
\end{equation}
where
\begin{align}
&\begin{cases}
a_{2m}=(2m-\alpha-\beta)(2m-\alpha)(2m-\beta)/4,\\
a_{2m-1}=(2m-2\!+\!\alpha\!+\!\beta)(2m-2\!+\!\alpha)(2m-2\!+\!\beta)/4,
\end{cases}\\
&\begin{cases}
b_{2m}(x)=x\!+\!2(1\!-\!\alpha\!-\!\beta),\\
b_{2m-1}(x)=(2m-1)x.
\end{cases}
\end{align}
Throughout this section, we use $\omega$ to denote $x+2(1-\alpha-\beta)$. We now assume that the adjoint factors $\{\phi_m\}_{m=1}^{\infty}$ are defined by
\begin{equation}\begin{cases}
\phi_{2m}=a_{2m}-r_{2m-1}\left(b_{2m}(x)+r_{2m}\right),\\
\phi_{2m-1}=a_{2m-1}-r_{2m-2}\left(b_{2m-1}(x)+r_{2m-1}\right),
\end{cases}\end{equation}
where modifying factors have the form of $r_{2m}=u_1 m^2+v_1 m+w_1, r_{2m-1}=u_2 m^2+v_2 m+w_2$.

In order to prove Theorem 3 by the Bauer-Muir transformation, the six parameters\\
$(u_1,u_2,v_1,v_2,w_1,w_2)$ will be so chosen such that both $\phi_{2m}$ and $\phi_{2m-1}$ are always constants, respectively, i.e. $[m^j]\phi_{2m}=0$ and $[m^j]\phi_{2m-1}=0$ for $j=1,2,3,4$. 
With the aid of \emph{Mathematica}, we choice
\begin{equation}
~~\left(u_1,u_2,v_1,v_2,w_1,w_2\right)=\Bigl(0,-\frac{2 x}{\omega},-\frac{\omega}{x},-x+\frac{2x}{\omega},-\frac{\omega}{2},-\frac{xh}{2\omega} \Bigr),~~
\end{equation}
where $h=(x-\alpha-\beta)^2+\alpha\beta$. In this case, it is not difficult to check that the following relations hold
\begin{equation}
\phi_{2m}=\frac 14(x-\alpha)(x-\alpha-\beta)(x-\beta)=-\phi_{2m-1}\ne 0.
\end{equation}
We evidently observe that the identity $\phi_{m+1}/\phi_{m}=-1$ always holds for any $m\in\mathbb{N}$. In addition, $b_0(x)/2+r_0=0,\quad b_1(x)+r_1=w_2,$ and
\begin{align}\begin{cases}
b_{2m}(x)+r_{2m}+r_{2m-2}=-(2m-1)\frac{\omega}{x},\\ b_{2m+1}(x)+r_{2m+1}+r_{2m-1}=-(4m^2+h)\frac{x}{\omega}.
\end{cases}\end{align}
It follows from the Bauer-Muir transformation and an equivalence transformation~(see~\cite[p.~73, Theorem 9]{LW}) that
\begin{align}
F(x)=&\frac{b_0(x)}{2}+r_0+\label{Theorem 3-RL-1}\\
&\begin{array}{cccccccc}
\phi_1& & -a_1& & -a_2 &\\
\cline{1-1}\cline{3-3}\cline{5-5}\cline{7-7}
b_1(x)+ r_1&\!+\!&b_2(x)+r_2+r_0&\!+\!&b_3(x)+r_3+r_1&\!+\cdots
\end{array}\nonumber\\
=&\begin{array}{ccccccc}
\phi_1& &-a_1& &-a_2 &\\
\cline{1-1}\cline{3-3}\cline{5-5}
w_2&\!+\!&-1\cdot\omega/x&\!+\!&-(4+h)\omega/x&\!+\!\cdots\!+\end{array}\nonumber\\
&\begin{array}{cccccc}
-a_{2m-1}&& -a_{2m}&\\
\cline{1-1}\cline{3-3}\cline{5-5}
-(2m-1)\omega/x&\!+\!&-(4m^2+h)\omega/x&\!+\cdots
\end{array}\nonumber\\
=&\begin{array}{ccccccc}
\phi_1\omega/x& & a_1& & a_2 &\\
\cline{1-1}\cline{3-3}\cline{5-5}
 h/2&\!+\!&1&\!+\!&-(4+h)&\!+\!\cdots\!+
\end{array}\nonumber\\
&\begin{array}{cccccc}
a_{2m-1}&& a_{2m}&\\
\cline{1-1}\cline{3-3}\cline{5-5}
2m-1&\!+\!&-(4m^2+h)&\!+\cdots
\end{array}\nonumber\\
:=&\frac{\phi_1\omega/x}{h/2+}\K_{m=1}^{\infty}\Bigl(\frac{a_m}{d_m(x)}\Bigr),\quad (\mbox{write}).\nonumber
\end{align}
Secondly, to treat the continued fraction $h/2+\k_{m=1}^{\infty}\bigl(a_m/d_m(x)\bigr)$, we shall employ the Bauer-Muir transformation again. To do that, we let
\begin{equation}\begin{cases}
\psi_{2m}=a_{2m}-\gamma_{2m-1}\left(d_{2m}(x)+\gamma_{2m}\right),\\
\psi_{2m-1}=a_{2m-1}-\gamma_{2m-2}\left(d_{2m-1}(x)+\gamma_{2m-1}\right),
\end{cases}
\end{equation}
where $\gamma_{2m}=p_1 m^2+q_1 m+t_1, \gamma_{2m-1}=p_2 m^2+q_2 m+t_2$. On taking
\begin{equation}
(p_1,p_2,q_1,q_2,t_1,t_2)=\Bigl(2,0,x,-1,\frac h2,-\frac x2+\alpha+\beta\Bigr),
\end{equation}
we can show that
\begin{align}\begin{cases}
\psi_{2m}=-\frac 14(x-\alpha-2\beta)(x-\alpha-\beta)(x-2\alpha-\beta)=-\psi_{2m-1}\ne 0,\\
d_{2m}(x)+\gamma_{2m}+\gamma_{2m-2}=(2m-1)(x-2),\\
d_{2m+1}(x)+\gamma_{2m+1}+\gamma_{2m-1}=-(\omega-2).
\end{cases}\end{align}
We find easily that the identity $\psi_{m+1}/\psi_{m}=-1$ always holds for any $m\in\mathbb{N}$. By using the Bauer-Muir transformation and an equivalence transformation again, one deduces that
\begin{align}
&\frac h2+\K_{m=1}^{\infty}\Bigl(\frac{a_m}{d_m(x)}\Bigr)\label{Theorem 3-RL-2}\\
=&\frac h2+\gamma_0+\begin{array}{cccccccc}
\psi_1& & -a_1& & -a_2 &\\
\cline{1-1}\cline{3-3}\cline{5-5}\cline{7-7}
d_1(x)+ \gamma_1&\!+\!&d_2(x)+\gamma_2+\gamma_0&\!+\!&d_3(x)+\gamma_3+\gamma_1&\!+\cdots
\end{array}\nonumber\\
=&\begin{array}{ccccccc}
\psi_1& &-a_1& &-a_2 &\\
\cline{1-1}\cline{3-3}\cline{5-5}
-(\omega-2)/2&\!+\!&1\cdot(x-2)&\!+\!&-(\omega-2)&\!+\!\cdots\!+
\end{array}\nonumber\\
&\begin{array}{cccccc}
-a_{2m-1}&& -a_{2m}&\\
\cline{1-1}\cline{3-3}\cline{5-5}
(2m-1)(x-2)&\!+\!&-(\omega-2)&\!+\cdots\!
\end{array}\nonumber\\
=&\begin{array}{ccccccc}
-\psi_1& &a_1& &a_2 &\\
\cline{1-1}\cline{3-3}\cline{5-5}
(\omega-2)/2&\!+\!&1\cdot(x-2)&\!+\!&\omega-2&\!+\!\cdots\!+
\end{array}\nonumber\\
&\begin{array}{cccccc}
a_{2m-1}&& a_{2m}&\\
\cline{1-1}\cline{3-3}\cline{5-5}
(2m-1)(x-2)&\!+\!&\omega-2&\!+\cdots
\end{array}\nonumber\\
=&\frac{-\psi_1}{F(x-2)}.\nonumber
\end{align}
Thirdly, substituting \eqref{Theorem 3-RL-2} into \eqref{Theorem 3-RL-1}, then replacing $x$ by $x+2$, after some simplifications we obtain that
\begin{equation}
\frac{F(x)}{F(x+2)}=\frac{\left(x+2-\alpha-2\beta\right)(x+2)\left(x+2-2\alpha-\beta\right)}
{\left(x+2-\alpha\right)\left(x+4-2\alpha-2\beta\right)\left(x+2-\beta\right)}.
\end{equation}
By iteration of this formula, we derive that, for each $m\in\mathbb{N}$
\begin{align}
&\frac{F(x)}{F(x+2m)}\\
=&\frac{\prod_{k=0}^{m-1}\left(x+2k+2-\alpha-2\beta\right)\left(x+2k+2\right)\left(x+2k+2-2\alpha-\beta\right)}
{\prod_{k=0}^{m-1}\left(x+2k+2-\alpha\right)\left(x+2k+4-2\alpha-2\beta\right)\left(x+2k+2-\beta\right)}\nonumber\\
=&\frac{\prod_{k=0}^{m-1}\left(\frac 12(x+2-\alpha-2\beta)+k\right)\left(\frac 12 (x+2)+k\right)\left(\frac 12(x+2-2\alpha-\beta)+k\right)}
{\prod_{k=0}^{m-1}\left(\frac 12(x+2-\alpha)+k\right)\left(\frac 12(x+4-2\alpha-2\beta)+k\right)\left(\frac 12(x+2-\beta)+k\right)}\nonumber\\
=&\frac {1}{m} \frac{\prod_{k=0}^{m-1}\left(\frac 12(x+2-\alpha-2\beta)+k\right)(m!)^{-1}m^{-\frac 12(x+2-\alpha-2\beta)}}{\prod_{k=0}^{m-1}\left(\frac 12(x+2-\alpha)+k\right)(m!)^{-1}m^{-\frac 12(x+2-\alpha)}}\times\nonumber\\
& \frac{\prod_{k=0}^{m-1}\left(\frac 12(x+2)+k\right)(m!)^{-1}m^{-\frac 12(x+2)}}{\prod_{k=0}^{m-1}\left(\frac 12(x+4-2\alpha-2\beta)+k\right)(m!)^{-1}m^{-\frac 12(x+4-2\alpha-2\beta)}}\times\nonumber\\
& \frac{\prod_{k=0}^{m-1}\left(\frac 12(x+2-2\alpha-\beta)+k\right)(m!)^{-1}m^{-\frac 12(x+2-2\alpha-\beta)}}{\prod_{k=0}^{m-1}\left(\frac 12(x+2-\beta)+k\right)(m!)^{-1}m^{-\frac 12(x+2-\beta)}}.\nonumber
\end{align}
It follows from Euler's formula \eqref{Euler's formula} that
\begin{align}
&\lim_{m\rightarrow\infty}\frac{F(x)m}{F(x+2m)}\label{Theorem 3-3}\\
=&\frac{\Gamma\left(\frac 12 ( x+2 - \alpha)\right) \Gamma\left(\frac 12( x+4-2(\alpha+\beta))\right)\Gamma
 \left( \frac 12 (x+2 -\beta)\right)}{\Gamma\left(\frac 12 (x+2 -\alpha-2\beta)\right) \Gamma\left(\frac 12(x+2)\right)\Gamma\left(\frac 12 ( x +2- 2\alpha-\beta)\right)}.\nonumber
\end{align}
From the definition of $F(x)$ in \eqref{Theorem 3-F-D}, we note that
\begin{equation}
\lim_{m\rightarrow\infty}\frac{m}{F(x+2m)}=1.\label{Theorem 3-4}
\end{equation}
Finally, combining \eqref{Theorem 3-F-D}, \eqref{Theorem 3-3}, \eqref{Theorem 3-4}, and applying an equivalence transformation, we get \eqref{Theorem 3} at once.\qed

\section{The proof of Theorem 4}
\proof  We shall prove Theorem 4 for $0<\alpha,\beta<1$ and $x>2$. An argument of Jacobsen~\cite[p.~419--420, Theorem 2.3(iv)]{Jac} can then be employed to extend the domains of convergence for $x, l, n,$ and $\eta$ to those indicated.

Firstly, let us consider the following continued fraction
\begin{equation}
F(x)=\frac{b_0(x)}{2}+\K_{m=1}^{\infty}\Bigl(\frac{a_m}{b_m(x)}\Bigr),\label{Theorem 4-F-D}
\end{equation}
where
\begin{align}
&\begin{cases}
a_{2m}=\frac{(2m-\alpha-\beta)(2m-\alpha)(2m-\beta)}{4},\\
a_{2m-1}=\frac{(2m-2+\alpha+\beta)(2m-2+\alpha)(2m-2+\beta)}{4},
\end{cases}\\
&\begin{cases}
b_{2m}(x)=x,\\
b_{2m-1}(x)=(2m-1)\left(x+2(1-\alpha-\beta)\right).
\end{cases}
\end{align}
Throughout this section, we let $\omega=x+2(1-\alpha-\beta)$. Assume that the adjoint factors $\{\phi_m\}_{m=1}^{\infty}$ are defined as
\begin{equation}\begin{cases}
\phi_{2m}=a_{2m}-r_{2m-1}\left(b_{2m}(x)+r_{2m}\right),\\
\phi_{2m-1}=a_{2m-1}-r_{2m-2}\left(b_{2m-1}(x)+r_{2m-1}\right),
\end{cases}\end{equation}
where modifying factors have the form of $r_{2m}=u_1 m^2+v_1 m+w_1, r_{2m-1}=u_2 m^2+v_2 m+w_2$.

In order to prove Theorem 4, the six parameters $(u_1,u_2,v_1,v_2,w_1,w_2)$ will be so chosen such that both $\phi_{2m}$ and $\phi_{2m-1}$ are always constants, respectively, i.e. $[m^j]\phi_{2m}=0$ and $[m^j]\phi_{2m-1}=0$ for $j=1,2,3,4$. 
With the aid of \emph{Mathematica}, we take
\begin{equation}
~~\left(u_1,u_2,v_1,v_2,w_1,w_2\right)=\Bigl(0,\frac{2\omega}{x},\frac{x}{\omega},-\omega-\frac{(2+x)\omega}{x},-\frac x2,\frac{h\omega}{2x} \Bigr),~~
\end{equation}
where $h=(x-\alpha+2)^2+(-4-2x+3\alpha)\beta+\beta^2$. In this case, it is not difficult to check that the following relations hold
\begin{equation}
\phi_{2m}=-\frac 14(x+2-\alpha-2\beta)(x+2-\alpha-\beta)(x+2-2\alpha-\beta)=-\phi_{2m-1}\ne 0.
\end{equation}
We evidently observe that the identity $\phi_{m+1}/\phi_{m}=-1$ always holds for any $m\in\mathbb{N}$. Moreover, $b_0(x)/2+r_0=0,\quad b_1(x)+r_1=w_2,$ and
\begin{align}\begin{cases}
b_{2m}(x)+r_{2m}+r_{2m-2}=(2m-1)\frac{x}{\omega},\\
 b_{2m+1}(x)+r_{2m+1}+r_{2m-1}=(4m^2+h)\frac{\omega}{x}.
\end{cases}\end{align}
It follows from the Bauer-Muir transformation and an equivalence transformation~(see~\cite[p.~73, Theorem 9]{LW}) that
\begin{align}
F(x)=&\frac{b_0(x)}{2}+r_0+\label{Theorem 4-RL-1}\\
&\begin{array}{cccccccc}
\phi_1& & -a_1& & -a_2 &\\
\cline{1-1}\cline{3-3}\cline{5-5}\cline{7-7}
b_1(x)+ r_1&\!+\!&b_2(x)+r_2+r_0&\!+\!&b_3(x)+r_3+r_1&\!+\cdots
\end{array}\nonumber\\
=&\begin{array}{ccccccc}
\phi_1& &-a_1& &-a_2 &\\
\cline{1-1}\cline{3-3}\cline{5-5}
w_2&\!+\!&1\cdot\frac{x}{\omega}&\!+\!&(4+h)\frac{\omega}{x}&\!+\!\cdots\!+
\end{array}\nonumber\\
&\begin{array}{cccccc}
-a_{2m-1}&& -a_{2m}&\\
\cline{1-1}\cline{3-3}\cline{5-5}
(2m-1)\frac{x}{\omega}&\!+\!&(4m^2+h)\frac{\omega}{x}&\!+\cdots
\end{array}\nonumber\\
=&\begin{array}{ccccccc}
\phi_1\frac{x}{\omega}& & a_1& & a_2 &\\
\cline{1-1}\cline{3-3}\cline{5-5}
 h/2&\!+\!&-1&\!+\!&4+h&\!+\!\cdots\!+
\end{array}\nonumber\\
&\begin{array}{cccccc}
a_{2m-1}&& a_{2m}&\\
\cline{1-1}\cline{3-3}\cline{5-5}
-(2m-1)&\!+\!&4m^2+h&\!+\cdots
\end{array}\nonumber\\
:=&\frac{\phi_1\frac{x}{\omega}}{ h/2+}\K_{m=1}^{\infty}\Bigl(\frac{a_m}{d_m(x)}\Bigr),\quad (\mbox{write}).\nonumber
\end{align}
Secondly, to treat the continued fraction $h/2+\k_{m=1}^{\infty}\bigl(a_m/d_m(x)\bigr)$, we need to employ the Bauer-Muir transformation again. To do that, we let
\begin{equation}\begin{cases}
\psi_{2m}=a_{2m}-\gamma_{2m-1}\left(d_{2m}(x)+\gamma_{2m}\right),\\
\psi_{2m-1}=a_{2m-1}-\gamma_{2m-2}\left(d_{2m-1}(x)+\gamma_{2m-1}\right),
\end{cases}\end{equation}
where $\gamma_{2m}=p_1 m^2+q_1 m+t_1, \gamma_{2m-1}=p_2 m^2+q_2 m+t_2$. On taking
\begin{equation}
(p_1,p_2,q_1,q_2,t_1,t_2)=\Bigl(-2,0,\omega,1,-\frac h2,-1-\frac x2\Bigr),
\end{equation}
we have
\begin{align}\begin{cases}
\psi_{2m}=\frac 14(x+2-\alpha)(x+2-\alpha-\beta)(x+2-\beta)=-\psi_{2m-1}\ne 0,\\
d_{2m}+\gamma_{2m}+\gamma_{2m-2}=(2m-1)(2+\omega),\\
d_{2m+1}+\gamma_{2m+1}+\gamma_{2m-1}=-(x+2).
\end{cases}\end{align}
We find easily that the identity $\psi_{m+1}/\psi_{m}=-1$ always holds for any $m\in\mathbb{N}$. Applying the Bauer-Muir transformation and an equivalence transformation again, one derives that
\begin{align}
&\frac h2+\K_{m=1}^{\infty}\Bigl(\frac{a_m}{d_m(x)}\Bigr)\label{Theorem 4-RL-2}\\
=&\frac h2+\gamma_0+\begin{array}{cccccccc}
\psi_1& & -a_1& & -a_2 &\\
\cline{1-1}\cline{3-3}\cline{5-5}\cline{7-7}
d_1(x)+ \gamma_1&\!+\!&d_2(x)+\gamma_2+\gamma_0&\!+\!&d_3(x)+\gamma_3+\gamma_1&\!+\cdots
\end{array}\nonumber\\
=&\begin{array}{ccccccc}
\psi_1& &-a_1& &-a_2 &\\
\cline{1-1}\cline{3-3}\cline{5-5}
-(x+2)/2&\!+\!&1\cdot(2+\omega)&\!+\!&-(x+2)&\!+\!\cdots\!+
\end{array}\nonumber\\
&\begin{array}{cccccc}
-a_{2m-1}&& -a_{2m}&\\
\cline{1-1}\cline{3-3}\cline{5-5}
(2m-1)(2+\omega)&\!+\!&-(x+2)&\!+\cdots
\end{array}\nonumber\\
=&\begin{array}{ccccccc}
-\psi_1& &a_1& &a_2 &\\
\cline{1-1}\cline{3-3}\cline{5-5}
(x+2)/2&\!+\!&1\cdot(2+\omega)&\!+\!&x+2&\!+\!\cdots\!+
\end{array}\nonumber\\
&\begin{array}{cccccc}
a_{2m-1}&& a_{2m}&\\
\cline{1-1}\cline{3-3}\cline{5-5}
(2m-1)(2+\omega)&\!+\!&x+2&\!+\cdots
\end{array}\nonumber\\
=&\frac{-\psi_1}{F(x+2)}.\nonumber
\end{align}
Thirdly, substituting \eqref{Theorem 4-RL-2} into \eqref{Theorem 4-RL-1}, then after some simplifications we get
\begin{equation}
\frac{F(x)}{F(x+2)}=\frac{\left(x+2-\alpha-2\beta\right)x\left(x+2-2\alpha-\beta\right)}
{\left(x+2-\alpha\right)\left(x+2-2\alpha-2\beta\right)\left(x+2-\beta\right)}.
\end{equation}
By iteration of the foregoing identity, we obtain that, for each $m\in\mathbb{N}$
\begin{align}
&\frac{F(x)}{F(x+2m)}\\
=&\frac{\prod_{k=0}^{m-1}\left(x+2k+2-\alpha-2\beta\right)\left(x+2k\right)\left(x+2k+2-2\alpha-\beta\right)}
{\prod_{k=0}^{m-1}\left(x+2k+2-\alpha\right)\left(x+2k+2-2\alpha-2\beta\right)\left(x+2k+2-\beta\right)}\nonumber\\
=&\frac{\prod_{k=0}^{m-1}\left(\frac 12(x+2-\alpha-2\beta)+k\right)\left(\frac 12 x+k\right)\left(\frac 12(x+2-2\alpha-\beta)+k\right)}
{\prod_{k=0}^{m-1}\left(\frac 12(x+2-\alpha)+k\right)\left(\frac 12(x+2-2\alpha-2\beta)+k\right)\left(\frac 12(x+2-\beta)+k\right)}\nonumber\\
=&\frac {1}{m} \frac{\prod_{k=0}^{m-1}\left(\frac 12(x+2-\alpha-2\beta)+k\right)(m!)^{-1}m^{-\frac 12(x+2-\alpha-2\beta)}}{\prod_{k=0}^{m-1}\left(\frac 12(x+2-\alpha)+k\right)(m!)^{-1}m^{-\frac 12(x+2-\alpha)}}\times\nonumber\\
&\frac{\prod_{k=0}^{m-1}\left(\frac 12x+k\right)(m!)^{-1}m^{-\frac 12x}}{\prod_{k=0}^{m-1}\left(\frac 12(x+2-2\alpha-2\beta)+k\right)(m!)^{-1}m^{-\frac 12(x+2-2\alpha-2\beta)}}\times\nonumber\\
&\frac{\prod_{k=0}^{m-1}\left(\frac 12(x+2-2\alpha-\beta)+k\right)(m!)^{-1}m^{-\frac 12(x+2-2\alpha-\beta)}}{\prod_{k=0}^{m-1}\left(\frac 12(x+2-\beta)+k\right)(m!)^{-1}m^{-\frac 12(x+2-\beta)}}.\nonumber
\end{align}
Applying Euler's formula \eqref{Euler's formula}, we find
\begin{align}
&\lim_{m\rightarrow\infty}\frac{F(x)m}{F(x+2m)}\label{Theorem 4-3}\\
=&\frac{\Gamma\left(\frac 12 ( x+2 - \alpha)\right) \Gamma\left(\frac 12( x+2-2(\alpha+\beta))\right)\Gamma
 \left( \frac 12 (x+2 -\beta)\right)}{ \Gamma\left(\frac 12 ( x+2 - 2\alpha-\beta)\right) \Gamma\left(\frac 12x\right)\Gamma
  \left(\frac 12 (x+2 -\alpha-2\beta)\right)}.\nonumber
\end{align}
From the definition of $F(x)$ in \eqref{Theorem 4-F-D}, we easily show that
\begin{equation}
\lim_{m\rightarrow\infty}\frac{m}{F(x+2m)}=1.\label{Theorem 4-4}
\end{equation}
Finally, combining \eqref{Theorem 4-F-D}, \eqref{Theorem 4-3} and \eqref{Theorem 4-4}, then using an equivalence transformation, and this will complete the proof of Theorem 4.\qed
\section{The proof of Theorem 5}

We shall prove Theorem 5 for $\alpha,\beta>0$, $\alpha+\beta<2$,  and $x>2$. The parabola theorem of Jacobsen~\cite[p.~419--420, Theorem 2.3(v)]{Jac} can then be employed to extend the domains of convergence for $x, \alpha$ and $\beta$ to those indicated. Throughout this section, for brevity we write $\sigma=\frac 14\alpha\beta(\alpha+\beta)$.

Firstly, let us consider the following generalized continued fraction
\begin{equation}
F(x):=b_0(x)+\K_{m=1}^{\infty}\Bigl(\frac{a_m}{b_m(x)}\Bigr),\label{conjecture 1-F-def}
\end{equation}
where two sequences $\{a_m\}_{m=1}^{\infty}$ and $\{b_m(x)\}_{m=0}^{\infty}$ are given as Theorem 5.
To employ the Bauer-Muir transformation, we assume that \emph{modifying factors} $\{r_m\}_{m=0}^{\infty}$ are of the form
\begin{equation}\begin{cases}
r_{2m}=u_1 m^2+v_1 m+w_1,\\
r_{2m-1}=u_2 m^2+v_2 m+w_2,
\end{cases}\end{equation}
 to be specified below. Let the \emph{adjoint factors} $\{\Phi_m\}_{m=1}^{\infty}$ for the continued fraction $F(x)$ be defined as follows
\begin{equation}\begin{cases}
\Phi_{2m}=a_{2m}-r_{2m-1}\left(b_{2m}(x)+r_{2m}\right),\\
\Phi_{2m-1}=a_{2m-1}-r_{2m-2}\left(b_{2m-1}(x)+r_{2m-1}\right).
\end{cases}\end{equation}
In order to prove Theorem 5, the six parameters $(u_1,u_2,v_1,v_2,w_1,w_2)$ will be so chosen such that both $\Phi_{2m}$ and $\Phi_{2m-1}$ are always constants, respectively, i.e. $[m^j]\Phi_{2m}=0$ and $[m^j]\Phi_{2m-1}=0$ for $j=1,2,3,4$. 
With the aid of \emph{Mathematica}, we take 
\begin{align}
&\left(u_1,u_2,v_1,v_2,w_1,w_2\right)\label{conjecture 1-first solution}\\
=&\Bigl(2(1+x),0,3+2x-x^2,\frac {1}{1+x},\frac 12(1+x)(3+x^2-\alpha^2-\alpha\beta-\beta^2),-\frac 12\Bigr).
\nonumber
\end{align}
In this case, one may verify that the following relations hold
\begin{equation}
\Phi_{2m}=\frac 14(x+1+\alpha)(x+1-\alpha-\beta)(x+1+\beta)=-\Phi_{2m-1}\ne 0.
\end{equation}
We evidently observe that the identity $\Phi_{m+1}/\Phi_{m}=-1$ always holds for any $m\in\mathbb{N}$. Moreover, $b_0(x)+r_0=b_0(x)+w_1,~ b_1(x)+r_1=\frac 12+\frac {1}{x+1}$, and
\begin{align}\begin{cases}
b_{2m}(x)+r_{2m}+r_{2m-2}=(x+1)(4m^2+\rho),\\
 b_{2m+1}(x)+r_{2m+1}+r_{2m-1}=\frac{2m+1}{x+1},
\end{cases}\end{align}
where $\rho=(x+1)^2-\alpha^2-\alpha\beta-\beta^2$.
It follows from the Bauer-Muir transformation that
\begin{align*}
F(x)=&b_0(x)+r_0\\
&+\begin{array}{ccccccc}
\phi_1& & -a_1& & -a_k&\\
\cline{1-1}\cline{3-3}\cline{5-5}\cline{7-7}
b_1(x)+ r_1&\!+\!&b_2(x)+r_2+r_0&\!\!+\!\cdots\!+\!\!&b_{k+1}(x)+r_{k+1}+r_{k-1}&\!+\!\cdots\!
\end{array}.
\end{align*}
By an equivalence transformation, we obtain that
\begin{align}
&F(x)-(b_0(x)+w_1)\label{conjecture-F-second expression}\\
=&\begin{array}{ccccccc}
\Phi_1& & -a_1& & -a_2 &\nonumber\\
\cline{1-1}\cline{3-3}\cline{5-5}\cline{7-7}
\frac 12+\frac {2}{x+1}&\!+\!&(x+1)(4\cdot 1^2+\rho)&\!+\!&\frac{2\cdot 1+1}{x+1}&\!+\!\cdots\!+
\end{array}\nonumber\\
&\begin{array}{ccccc}
-a_{2m-1}&&-a_{2m}&\\
\cline{1-1}\cline{3-3}\cline{5-5}
(x+1)(4m^2+\rho) &\!+\!&\frac{2m+1}{x+1}&\!+\!\cdots
\end{array}\nonumber\\
=&\begin{array}{ccccccc}
(x+1)\Phi_1& &-a_1& &-a_2 &\\
\cline{1-1}\cline{3-3}\cline{5-5}
(x+3)/2&\!+\!&4\cdot 1^2+\rho&\!+\!&2\cdot 1+1&\!+\!\cdots\!+\!
\end{array}\nonumber\\
&\begin{array}{cccccc}
-a_{2m-1}&& -a_{2m}&\\
\cline{1-1}\cline{3-3}\cline{5-5}
4m^2+\rho&\!+\!&2m+1&\!+\cdots
\end{array}\nonumber\\
=&\begin{array}{ccccccccccc}
(x+1)\Phi_1& &a_1& &a_2 &\\
\cline{1-1}\cline{3-3}\cline{5-5}
(x+3)/2&\!+\!&-4\cdot 1^2-\rho&\!+\!&2\cdot 1+1&\!+\!\cdots\!+
\end{array}\nonumber\\
&\begin{array}{cccccc}
a_{2m-1}&& a_{2m}&\\
\cline{1-1}\cline{3-3}\cline{5-5}
-4m^2-\rho&\!+\!&2m+1&\!+\cdots
\end{array}\nonumber\\
:=&\frac{(x+1)\Phi_1}{(x+3)/2+}\K_{m=1}^{\infty}\Bigl(\frac{a_m}{d_m(x)}\Bigr),\quad (\mbox{write}).\nonumber
\end{align}
Secondly, to treat the continued fraction $(x+3)/2+\k_{m=1}^{\infty}\bigl(a_m/d_m(x)\bigr)$, we need to apply the Bauer-Muir transformation again. To this end, we let
\begin{equation}\begin{cases}
\Psi_{2m}=a_{2m}-\gamma_{2m-1}\left(d_{2m}(x)+\gamma_{2m}\right),\\
\Psi_{2m-1}=a_{2m-1}-\gamma_{2m-2}\left(d_{2m-1}(x)+\gamma_{2m-1}\right),
\end{cases}\end{equation}
where $\gamma_{2m}=p_1 m^2+q_1 m+t_1, \gamma_{2m-1}=p_2 m^2+q_2 m+t_2$. In the same way as the above approach, on taking
\begin{equation}
(p_1,p_2,q_1,q_2,t_1,t_2)=\Bigl(0,2,-1,-x-1,\frac {x-1}{2},\frac{(x+1)^2-\alpha^2\alpha\beta-\beta^2}{2}\Bigr),\label{conjecture 1-second solution}
\end{equation}
we find that
\begin{align}\begin{cases}
\Psi_{2m}=-\frac 14(x+1-\alpha)(x+1+\alpha+\beta)(x+1-\beta)=-\Psi_{2m-1}\ne 0,\\
d_{2m}(x)+\gamma_{2m}+\gamma_{2m-2}=x+1,\\
d_{2m+1}(x)+\gamma_{2m+1}+\gamma_{2m-1}=-(x+3)(2m+1),\\
\frac{x+3}{2}+\gamma_0=x+1,\quad d_1(x)+\gamma_1=-\frac 12\big(7+4x+x^2-\alpha^2-\alpha\beta-\beta^2\big).
\end{cases}
\end{align}
We see easily that the identity $\Psi_{m+1}/\Psi_{m}=-1$ always holds for any $m\in\mathbb{N}$. Using the Bauer-Muir transformation again, and then employing an equivalence transformation, one derives that
\begin{align}
&\frac{x+3}{2}+\K_{m=1}^{\infty}\Bigl(\frac{a_m}{d_m(x)}\Bigr)\label{conjecture-GG-second expression}\\
=&\frac {x+3}{2}+\gamma_0+\begin{array}{ccccccc}
\Psi_1& & -a_1& & -a_2 &\\
\cline{1-1}\cline{3-3}\cline{5-5}\cline{7-7}
d_1(x)+\gamma_1&\!+\!&d_2(x)+\gamma_2+\gamma_0&\!+\!&d_3(x)+\gamma_3+\gamma_1&\!+\!\cdots\!
\end{array}\nonumber\\
=&x+1+\begin{array}{ccccccc}
\Psi_1& & -a_1& & -a_2 &\nonumber\\
\cline{1-1}\cline{3-3}\cline{5-5}\cline{7-7}
d_1(x)+\gamma_1&\!+\!&x+1&\!+\!&-(x+3)(2\cdot 1+1)&\!+\!\cdots\!+\!
\end{array}\nonumber\\
&\begin{array}{ccccc}
-a_{2m-1}&&-a_{2m}&\\
\cline{1-1}\cline{3-3}\cline{5-5}
x+1 &\!+\!&-(x+3)(2m+1)&\!+\!\cdots\!
\end{array}\nonumber\\
=&x+1+\begin{array}{ccccccc}
-\Psi_1& & a_1& & a_2 &\nonumber\\
\cline{1-1}\cline{3-3}\cline{5-5}\cline{7-7}
-(d_1(x)+\gamma_1)+&\!+\!&x+1&\!+\!&(x+3)(2\cdot 1+1)&\!+\!\cdots\!+\!
\end{array}\nonumber\\
&\begin{array}{ccccc}
a_{2m-1}&&a_{2m}&\\
\cline{1-1}\cline{3-3}\cline{5-5}
x+1 &\!+\!&(x+3)(2m+1)&\!+\!\cdots\!
\end{array}\nonumber\\
=&x+1+\begin{array}{ccccccc}
-(x+1)\Psi_1& & a_1& & a_2 \nonumber\\
\cline{1-1}\cline{3-3}\cline{5-5}\cline{7-7}
-(x+1)(d_1(x)+\gamma_1)+&\!+\!&1&\!+\!&(x+1)(x+3)(2\cdot 1+1)
\end{array}\nonumber\\
&\begin{array}{ccccccc}
&a_{2m-1}&&a_{2m}&\\
\cline{2-2}\cline{4-4}\cline{6-6}
\!+\!\cdots\!+\!&1 &\!+\!&(x+1)(x+3)(2m+1)&\!+\!\cdots\!
\end{array}\nonumber\\
=&x+1+\frac{-(x+1)\Psi_1}{F(x+2)+b_0(x)+w_1+2\sigma},\nonumber
\end{align}
here in the last relation we used $(x+1)(x+3)=(x+2)^2-1$ and $-(x+1)(d_1(x)+\gamma_1)=b_0(x+2)+\left(b_0(x)+w_1+2\sigma\right)$.

Thirdly, substituting \eqref{conjecture-GG-second expression} into \eqref{conjecture-F-second expression}, after some simplifications we obtain
\begin{align}\label{recurrence relation of F(x)}
F(x)&= b_0(x)+w_1+\frac{(x+1)\Phi_1}{x+1+\frac{-(x+1)\Psi_1}{F(x+2)+b_0(x)+w_1+2\sigma}}\\
  &=b_0(x)+w_1+\frac{\Phi_1\left(F(x+2)+b_0(x)+w_1+2\sigma\right)}{F(x+2)+b_0(x)+w_1+2\sigma-\Psi_1}\nonumber\\
  &=b_0(x)+w_1+\frac{\Phi_1\left(F(x+2)+b_0(x)+w_1+2\sigma\right)}{F(x+2)+\frac 14\left((x+1)^2-\alpha^2-\alpha\beta-\beta^2\right)}.\nonumber
\end{align}
We now introduce a function $Q(x)$ to be defined as follows
\begin{equation}
  Q(x):=\frac{F(x)-\sigma}{F(x)+\sigma}\Leftrightarrow \frac{1-Q(x)}{1+Q(x)}=\frac{\sigma}{F(x)}.\label{definition of Q(x)}
\end{equation}
Via \emph{Mathematica}, it follows from \eqref{definition of Q(x)} and \eqref{recurrence relation of F(x)} that
\begin{align}
  Q(x) &=\frac{F(x)-\sigma}{F(x)+\sigma}\label{Recurrence relation for Q(x)}\\
  &=\frac{b_0(x)+w_1+\frac{\Phi_1\left(F(x+2)+b_0(x)+w_1+2\sigma\right)}{F(x+2)+\frac 14\left((x+1)^2-\alpha^2-\alpha\beta-\beta^2\right)}-\sigma}{b_0(x)+w_1+\frac{\Phi_1\left(F(x+2)+b_0(x)+w_1+2\sigma\right)}{F(x+2)+\frac 14\left((x+1)^2-\alpha^2-\alpha\beta-\beta^2\right)}+\sigma}\nonumber\\
  & =\frac{(x+1+\alpha)(x+1-\alpha-\beta)(x+1+\beta)}{(x+1-\alpha)(x+1+\alpha+\beta)(x+1-\beta)}\cdot\frac{F(x+2)-\sigma}{F(x+2)-\sigma}
  \nonumber\\
  &=\frac{(x+1+\alpha)(x+1-\alpha-\beta)(x+1+\beta)}{(x+1-\alpha)(x+1+\alpha+\beta)(x+1-\beta)}Q(x+2).\nonumber
\end{align}
By iteration of the foregoing relation, we find that, for each $m\in\mathbb{N}$
\begin{align}
&\frac{Q(x)}{Q(x+2m)}\\
=&\frac{\prod_{k=0}^{m-1}\left(x+2k+1+\alpha\right)\left(x+2k+1-\alpha-\beta\right)\left(x+2k+1+\beta\right)}
{\prod_{k=0}^{m-1}\left(x+2k+1-\alpha\right)(x+2k+1+\alpha+\beta)\left(x+2k+1-\beta\right)}\nonumber
\\
=&\frac{\prod_{k=0}^{m-1}\left(\frac 12(x+\alpha+1)+k\right)\left(\frac 12(x-\alpha-\beta+1)+k\right)\left(\frac 12(x+\beta+1)+k\right)}{\prod_{k=0}^{m-1}\left(\frac 12(x-\alpha+1)+k\right)\left(\frac 12(x+\alpha+\beta+1)+k\right)\left(\frac 12(x-\beta+1)+k\right)}\nonumber\\
=&\frac{\prod_{k=0}^{m-1}\left(\frac 12(x+\alpha+1)+k\right)(m!)^{-1}m^{-\frac 12(x+\alpha+q)}}{\prod_{k=0}^{m-1}\left(\frac 12(x-\alpha+1)+k\right)(m!)^{-1}m^{-\frac 12(x-\alpha+1)}}
\nonumber\\
&\times \frac{\prod_{k=0}^{m-1}\left(\frac 12(x-\alpha-\beta+1)+k\right)(m!)^{-1}m^{-\frac 12(x-\alpha-\beta+1)}}{\prod_{k=0}^{m-1}\left(\frac 12(x+\alpha+\beta+1)+k\right)(m!)^{-1}m^{-\frac 12(x+\alpha+\beta+1)}}\nonumber\\
&\times\frac{\prod_{k=0}^{m-1}\left(\frac 12(x+\beta+1)+k\right)(m!)^{-1}m^{-\frac 12(x+\beta+1)}}{\prod_{k=0}^{m-1}\left(\frac 12(x-\beta+1)+k\right)(m!)^{-1}m^{-\frac 12(x-\beta+1)}}.\nonumber
\end{align}
Note that $\lim_{m\rightarrow\infty}Q(x+2m)=1$ from the definition of $Q(x)$ in \eqref{definition of Q(x)} and \eqref{conjecture 1-F-def}. By means of Euler's formula~(see~\cite[p.~255, (6.1.2)]{AS})
\begin{equation}
\Gamma(z)=\lim_{n\rightarrow\infty}\frac{n!~n^z}{z(z+1)\cdots(z+n)}\quad(z\neq 0,-1,-2,\ldots),
\end{equation}
let $m$ tend to infinity, we derive that
\begin{equation}
Q(x)=\frac{ \Gamma\left(\frac 12 ( x - \alpha+1)\right) \Gamma\left(\frac 12(x+\alpha+\beta+1)\right)\Gamma
  \left(\frac 12 (x -\beta+1)\right)}{
\Gamma\left(\frac 12 ( x +\alpha+1)\right) \Gamma\left(\frac 12( x-\alpha-\beta+1)\right)\Gamma
 \left( \frac 12 (x +\beta+1)\right)}.\label{gamma expression for Q(x)}
\end{equation}
Combining \eqref{conjecture 1-F-def}, \eqref{definition of Q(x)} and \eqref{gamma expression for Q(x)}, then applying an equivalence transformation, we shall finish the proof of \eqref{first assertion of theorem 5}. Finally, by an equivalence transformation, we can rewrite the continued fraction on the right side of \eqref{first assertion of theorem 5} into the form of $\frac{\sigma}{b_0(x)+}\k_{m=1}^{\infty}\bigl(\frac{p_m}{q_m(x)}\bigr)$ with $q_{2m-1}(x)=1$ and $q_{2m}(x)=x^2-1$, since the continued fraction on right hand side of \eqref{second assertion of theorem 5} is the even part of the previous one, this completes the proof.\qed

\begin{remark}
Via the MC-algorithm, we can only guess that \eqref{second assertion of theorem 5} holds. However, we can not prove such a claim directly by the Bauer-Muir transformation, we have to find some ``simpler" continued fractions such that the continued fraction on the right side of \eqref{second assertion of theorem 5} is the even part of it. Therefore, firstly we should discuss 20 possibilities, and then pick out 3 ``uncomplicated" expressions among them. Secondly, write each of them into the polynomial continued fraction. Thirdly, try to check whether each of them satisfies the  identity \eqref{first assertion of theorem 5}. It is eventually shown by the Bauer-Muir transformation that only one among them is workable. The other two continued fractions can be are expressed as follows: for $j=1,2$
\begin{equation}
F_j(x)=b_0(1;x)+\K_{m=1}^{\infty}\Bigl(\frac{a_m(j}{b_m(j;x)}\Bigr),\label{conjecture 1-F1-def}
\end{equation}
where $b_0(1;x)=x^2-(1-\alpha)^2-\alpha\beta(\alpha+\beta)/4$, and for $k\in\mathbb{N}$
\begin{align}
&\begin{cases}
a_{2m}(1)=\frac{1}{4}(2m+\alpha)(2m+\alpha+\beta)(2m-\beta),\\
a_{2m-1}(1)=
\frac 14{(2m-\alpha)(2m-\alpha-\beta)(2m+\beta)},
\end{cases}\\
&\begin{cases}
b_{2m}(1;x)=(2m+1)\left(x^2-(1-\alpha)^2\right),\\
b_{2m-1}(1;x)=1;
 \end{cases}
\end{align}
and $b_0(2;x)=x^2-(1-\alpha-\beta)^2+\alpha\beta(\alpha+\beta)/4$, and for $k\in\mathbb{N}$
\begin{align}
&\begin{cases}
a_{2m}(2)=\frac{1}{4}(2m+\alpha)(2m+\alpha+\beta)(2m+\beta),\\
a_{2m-1}(2)=
\frac 14{(2m-\alpha)(2m-\alpha-\beta)(2m-\beta)},
\end{cases}\\
&\begin{cases}
b_{2m}(2;x)=(2m+1)\left(x^2-(1-\alpha-\beta)^2\right),\\
b_{2m-1}(2;x)=1,
\end{cases}
\end{align}
respectively.
\end{remark}
\section{Three open problems}


\begin{conjecture}
We let
\begin{equation}
P:=\frac{ \Gamma\left(\frac 12 ( x - \alpha+\frac 32)\right) \Gamma\left(\frac 12(x+\alpha+\beta+\frac 12)\right)\Gamma
  \left(\frac 12 (x -\beta+\frac 32)\right)}{
\Gamma\left(\frac 12 ( x +\alpha+\frac 32)\right) \Gamma\left(\frac 12( x-\alpha-\beta+\frac 12)\right)\Gamma
 \left( \frac 12 (x +\beta+\frac 32)\right)}.
\end{equation}
Suppose that $x$ is complex with $\Re(x)>0$, or assume that either one of $\alpha, \beta$ is an odd integer or $\alpha+\beta$ is an even integer. Then
\begin{equation}
\frac{1-P}{1+P}=\frac{(\alpha+\beta)/2}{x+\lambda_0+}\K_{m=1}^{\infty}\Bigl(\frac{\kappa_m}{x+\lambda_m}\Bigr),
\end{equation}
where
\begin{equation}
\begin{cases}\kappa_{2m}=\frac{(2m)^2-(\alpha+\beta)^2}{4},\\
\kappa_{2m-1}=\frac{\left((2m-1)^2-\alpha^2\right)\left((2m-1)^2-\beta^2\right)}{4(2m-1)^2},
\end{cases}
\quad \begin{cases}\lambda_{2m}=-\frac{\alpha\beta}{2(2m+1)},\\
\lambda_{2m-1}=\frac{\alpha\beta}{2(2m-1)}.
\end{cases}
\end{equation}
\end{conjecture}

\begin{conjecture}
Define
\begin{equation}
P:=\frac{ \Gamma\left(\frac 12 ( x - \alpha+\frac 12)\right) \Gamma\left(\frac 12(x+\alpha+\beta+\frac 32)\right)\Gamma
  \left(\frac 12 (x -\beta+\frac 12)\right)}{
\Gamma\left(\frac 12 ( x +\alpha+\frac 12)\right) \Gamma\left(\frac 12( x-\alpha-\beta+\frac 32)\right)\Gamma
 \left( \frac 12 (x +\beta+\frac 12)\right)}.
\end{equation}
Suppose that  $x$ is complex with $\Re(x)>0$, or assume that either one of $\alpha, \beta$ is an odd integer or $\alpha+\beta$ is an even integer. Then
\begin{equation}
\frac{1-P}{1+P}=-\frac{(\alpha+\beta)/2}{x+\lambda_0+}\K_{m=1}^{\infty}\Bigl(\frac{\kappa_m}{x+\lambda_m}\Bigr),
\end{equation}
where
\begin{equation}
\begin{cases}\kappa_{2m}=\frac{(2m)^2-(\alpha+\beta)^2}{4},\\
\kappa_{2m-1}=\frac{\left((2m-1)^2-\alpha^2\right)\left((2m-1)^2-\beta^2\right)}{4(2m-1)^2},
\end{cases}
\begin{cases}\lambda_{2m}=\frac{\alpha\beta}{2(2m+1)},\\
\lambda_{2m-1}=-\frac{\alpha\beta}{2(2m-1)}.
\end{cases}
\end{equation}

\end{conjecture}

\begin{conjecture} Let $x, l, n, \eta\in \mathbb{C}$, and $P$ be defined by
\begin{align}
&\frac{\Gamma\left(\frac 14(x+l+n+\eta+1)\right)\Gamma\left(\frac 14(x+l-n-\eta+1)\right)}
{\Gamma\left(\frac 14(x+l+n-\eta+3)\right)\Gamma\left(\frac 14(x+l-n+\eta+3)\right)}\times\\
&\frac{\Gamma\left(\frac 14(x-l+n-\eta+3)\right)\Gamma\left(\frac 14(x-l-n+\eta+3)\right)}{\Gamma
\left(\frac 14(x-l+n+\eta+1)\right)\Gamma\left(\frac 14(x-l-n-\eta+1)\right)}\nonumber.
\end{align}
Suppose that $x$ is complex with $\Re(x) >0$, or assume that either one of $n$ , $\eta$ is an odd integer or $l$ is an even integer. Then
\begin{equation}
\frac{1-P}{1+P}=\frac{l}{x-n\eta+}
\K_{m=1}^{\infty}\Bigl(\frac{a_m}{b_m(x)}\Bigr),\label{Conjecture 4}
\end{equation}
where
\begin{equation}\begin{cases}
a_{2m}=(2m)^2-l^2,\\ a_{2m-1}=\frac{\left((2m-1)^2-n^2\right)\left((2m-1)^2-\eta^2\right)}{(2m-1)^2},
\end{cases}
\begin{cases}
 b_{2m}(x)=x-\frac{n\eta}{2m+1},\\ b_{2m-1}(x)=x+\frac{n\eta}{2m-1}.
\end{cases}\end{equation}
\end{conjecture}
\begin{remark}
Note that there are four parameters in Conjecture 4. Setting $\eta=0$ , we obtain Entry 34 of Ramanujan on continued fractions involving gamma functions in his notebooks directly, see \cite[p.~156]{Berndt}. Applications for Theorem 5, Conjectures 1 to 3 will be discussed elsewhere.
\end{remark}
\begin{remark}
Of course, it is an interesting question to prove our theorems by the theory of hypergeometric functions or orthogonal polynomials~(for instance, see ~\cite{KLS,KS}). Moreover, it certainly would be a very interesting work to extend our results to $q$-analogue~(see \cite{GM1993,GM1998, Masson1995}). Finally, solving these conjectures should help progress the theory of continued fractions, the hypergeometric function, orthogonal polynomials, $q$-series, etc.
\end{remark}

\bibliographystyle{amsplain}

\end{document}